%% file: main.tex
\documentclass[11pt,a4paper]{article}
\input{preamble}

\title{Cost benefits of designing carbon capture processes\\ for part-load operation of power plants}
\author[1,2]{David Yang Shu}
\author[1]{Boxun Huang}
\author[1]{Yurim Kim}
\author[3]{Rahul Gandhi}
\author[2]{Randall Field}
\author[1]{Sungho Shin\thanks{Corresponding author: \texttt{sushin@mit.edu}.}}
\affil[1]{Department of Chemical Engineering, Massachusetts Institute of Technology}
\affil[2]{MIT Energy Initiative, Massachusetts Institute of Technology}
\affil[3]{Shell Technology Center}
\date{}

\begin{document}
\section*{Graphical abstract}
\begin{figure*}[!htb]
    \centering
    \includegraphics[width=15cm]{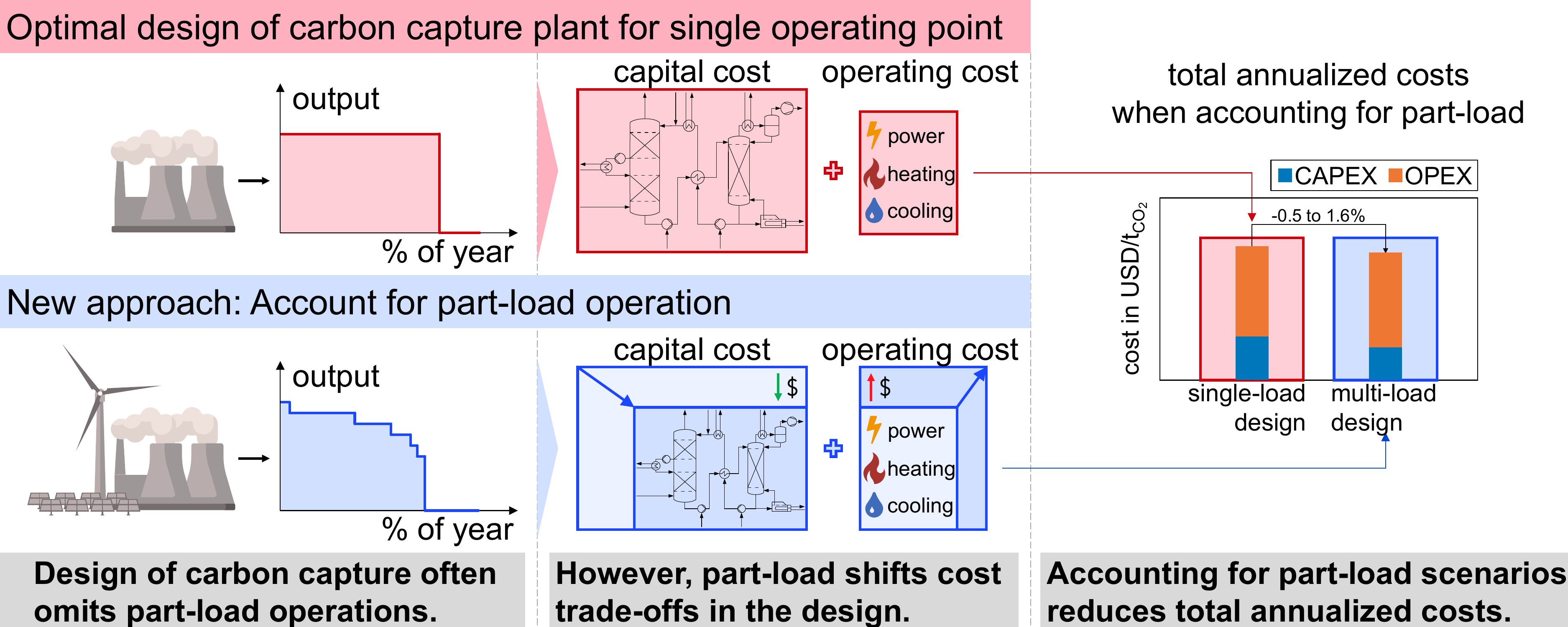}
    \end{figure*}
\newpage

\section*{Highlights}
\begin{itemize}
    \item Two-stage stochastic optimization of \ce{CO2} capture design with part-load scenarios
    \item Data-driven approach-to-equilibrium model simplifies model of carbon capture process
    \item Scenarios derived from real-world operating data of coal and NGCC power plants
    \item Multi-load designs lower capital cost by \SIrange{5}{9}{\percent} and annualized cost by \SIrange{0.5}{1.6}{\percent}
    \item Cost benefits are robust in sensitivity analyses and increase lower capacity factor
\end{itemize}

\maketitle

\begin{abstract}
As fossil power plants increasingly cycle to balance variable renewable generation, carbon capture plants must be designed for a wide range of part-load operating conditions, to avoid excessive capital costs. However, accounting for the variability in design optimization yields higher computational complexity, and is therefore not a current design methodology. Here, we apply two-stage stochastic optimization to identify optimal carbon capture plant designs accounting for part-load flue gas conditions. We reduce computational complexity through a data-driven approach-to-equilibrium model of the absorption and desorption processes. Operating scenarios are derived from real-world part-load operation data of coal and natural gas combined cycle power plants. Our results demonstrate that designs optimized over multiple operating scenarios substantially reduce equipment sizes, lowering capital costs by \SIrange{5}{9}{\percent} at the expense of higher operating costs, resulting in a moderate reduction in total annualized cost by \SIrange{0.5}{1.6}{\percent}. The lowest costs are achieved when targeting an average capture target across operating scenarios, which allows \ce{CO2} capture to be shifted from full-load to part-load periods to reduce peak equipment sizes without compromising overall \ce{CO2} emission reductions. These cost benefits are robust across sensitivity analyses and increase with higher shares of part-load operation. Thus, part-load variability should be considered at the design stage, particularly for plants subject to variable load.
\end{abstract}

\keywords{process systems engineering, stochastic programming, nonlinear programming, techno-economic assessment, carbon capture and storage (CCS)}

\section{Introduction}
\begin{figure}[!t]
    \centering
    \includegraphics[width=0.4\textwidth]{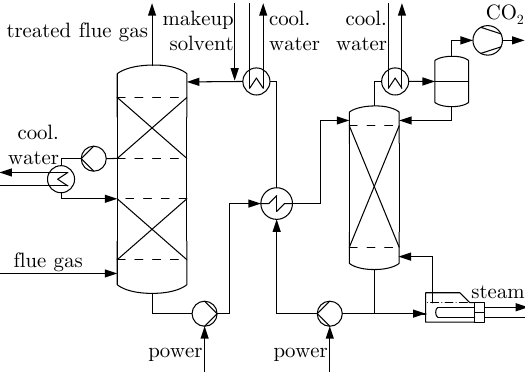}
    \caption{Process flow diagram of a solvent-based carbon capture process. }
    \label{fig:processflowdiagram}
\end{figure}

Fossil power plants still accounted for almost \SI{60}{\percent} of electricity generation in the United States in 2025 \citep{U.S.EnergyInformationAdministration.2026a}. A renewed interest in natural gas power plants to supply reliable power to data centers may lock in additional \ce{CO2} emissions over the lifetime of the plants. These \ce{CO2} emissions can be reduced by carbon capture processes. Mature capture processes today remove \ce{CO2} from flue gases through contact with a solvent in countercurrent flow in absorber columns (\Cref{fig:processflowdiagram}). The \ce{CO2} is removed from the rich solvent by heating it to produce a \ce{CO2}-rich stream that is further conditioned for permanent sequestration or utilization. 

Carbon capture processes are both capital- and energy intensive, resulting in high costs. Computational process system optimization identifies designs of carbon capture processes that balance capital and operating costs, while meeting process targets and technical and regulatory constraints. \citet{Chatziasteriou.2022} review works on simulation and optimization of carbon capture process designs. 

Such process designs should account for uncertainties and variabilities which are prevalent in operational conditions and in economic and technical parameters \citep{Li.2021}. One source of operational uncertainty in carbon capture processes in the variation in flue-gas conditions as power plant load changes. As fossil power plants provide grid-balancing flexibility, they are increasingly operated at part load in grids with expanding penetration of renewable power generation. Consequently, carbon capture plants should be designed for a wide range of part-load operating conditions, where they are operated at a fraction of their nameplate capacity.

However, accounting for the uncertainty resulting from variability in economic and technical parameters in process systems design can be computationally challenging \citep{Li.2021}. To reduce the complexity in carbon capture process design, rigorous handling of uncertainties is sometimes avoided by first obtaining deterministic designs from simulation \citep{Kuncheekanna.2025} or from optimization with detailed parametric studies \citep{Chatziasteriou.2022} and then applying design margins. Alternatively, computational complexity can be reduced by applying surrogate models instead of detailed process models \citep{Bahakim.2015, Giorgetti.2020}, or by focusing on the sizing of select components \citep{Cerrillo-Briones.2019, Zaman.2016}. Uncertainty can further be handled via robust optimization, yielding designs that guarantee feasibility for the considered uncertainty sets \citep{Cerrillo-Briones.2019, Sherman.2026, Akkor.2026} but may be overly conservative \citep{Li.2021}. Alternatively, uncertainty can be addressed by stochastic programming, which can be risk-neutral \citep{Li.2021}. Prior research considered the optimal sizing of the desorber and storage tanks via stochastic programming under uncertain electricity and \ce{CO2} prices or \ce{CO2} taxes \citep{Zaman.2016}. \citet{Seo.2023} further applied stochastic optimization to the design of carbon capture using three typical days as scenarios, focusing on dynamic operation to respond to a varying electricity price and demand. 

In contrast to previous studies, we apply two-stage stochastic optimization to isolate the impact of part-load flue gas variability as a design driver for capture plants.
We thus aim to answer the following question:

\begin{quote}
    \textit{How does accounting for the part-load operation prevalent in thermal power plants affect the optimal design and cost of solvent-based carbon capture processes?}
\end{quote}

To answer this question, we derive scenarios from real-world distributions of part-load operation levels for coal and \gls{ngcc} power plant units. We model the capture process as a set of interconnected unit operations in an equation-oriented formulation that enables us to simultaneously optimize the plant design and its operation in all scenarios. We include detailed techno-economic models that represent tradeoffs between performance, capital, and operating costs.

To keep the problem tractable, we describe the absorption and desorption processes in the columns using a data-drive approach-to-equilibrium model. This data-driven model is calibrated using simulation data to simplify the representation of rate limitations from mass transfer resistances and reaction kinetics. We further adopt a simplified thermodynamic model, assuming a constant heat of absorption of \ce{CO2} in \gls{mea}, and ideal-mixture enthalpy in both phases using pure-component models for heat capacity.

In the \Cref{sec:methods}, the process and cost models are described, before we present the case study assumptions for coal and \gls{ngcc} power plants in \Cref{sec:casestudy}. The results of the case studies including a sensitivity analysis are presented in \Cref{sec:results} and discussed in \ref{sec:discussions} before we conclude this paper (\Cref{sec:conclusions}).

\section{Methods}
\label{sec:methods}
We formulate the design problem of the post-combustion carbon capture system (\Cref{fig:processflowdiagram})  as a two-stage stochastic programming problem \cref{eq:problem} assuming a process based on \SI{30}{\wtpercent} aqueous \gls{mea} as solvent:
\begin{subequations}\label{eq:problem}
\begin{align}
\min_{\mathbf{x},\{\mathbf{y}_s\}_{s\in \mathcal{S}}} \quad
    & c^{\mathrm{toc}}(\mathbf{x}) + \sum_{s\in \mathcal{S}} w_s\, c^{\mathrm{op}}(\mathbf{x},\mathbf{y}_s) \label{eq:objective}\\
\text{s.t.} \quad
    & f_s(\mathbf{x},\mathbf{y}_s) \le 0, \, && s\in \mathcal{S},\label{eq:ineq}\\
    & h(\mathbf{x},\{\mathbf{y}_s\}_{s\in \mathcal{S}}) \le 0, \label{eq:ineq_coupling}\\
    & \mathbf{l_x} \le \mathbf{x} \le \mathbf{u_x},\quad \mathbf{l_y} \le \mathbf{y}_s \le \mathbf{u_y},\, && s\in \mathcal{S}.\label{eq:bounds}
\end{align}
\end{subequations}
Here, total annualized cost is chosen as the objective function \cref{eq:objective}, where $c^{\mathrm{toc}}(\mathbf{x})$ represents annualized total overnight cost, $c^{\mathrm{op}}(\mathbf{x},\mathbf{y}_s)$ the scenario-dependent annual operating costs, and $w_s$ is the weight for scenario $s\in\mathcal{S}$, where $\mathcal{S}$ is the set of operating scenarios. The first-stage design variables $\mathbf{x}$ include the equipment sizes of the capture process (\Cref{fig:processflowdiagram}), specifically the packing height and diameter of absorber and desorber columns, the maximum flow rates of pumps, the number of channels in the cross heat exchanger, the maximum compressor power, and the surface area of the heat exchangers for the reboiler, condenser, cooler, and intercooler. The number of channels in the cross heat exchanger is relaxed to a continuous variable to maintain a nonlinear programming formulation. 

The second-stage operating variables $\mathbf{y}_s$ are optimized to adapt to varying conditions in the scenarios $s$, and include the solvent loading and flow rate, and the reboiler duty. The scenarios differ in flue-gas composition and flow rate. The scenario-dependent constraints $f_s$ and the scenario-coupling constraints $h$ represent, e.g., techno-economic models, capture rate targets, and technical process models, such as thermodynamic properties, material and energy balances, or the connectivity of the unit operations in the carbon capture process. The first- and second-stage variables are bounded by $\mathbf{l_x}$, $\mathbf{u_x}$ and $\mathbf{l_y}$, $\mathbf{u_y}$, respectively. A detailed model description is provided in \Cref{SI:models}.
    
\subsection{Process system model}
The process system model determines the equipment sizing, utility demands, and captured amount of \ce{CO2}, and is composed of interconnected unit operations (\Cref{fig:processflowdiagram}). 

We adopt the following key assumptions in our process model: We assume ideal-mixture enthalpy in both phases, using the pure-component models for heat capacity, heat of vaporization and heat of absorption collected by \citet{Akula.2021}, neglecting heat of mixing.
For simplicity, we further neglect the solubility of \ce{O2} and \ce{N2} in the solvent and assume \gls{mea} to be nonvolatile. As an additional simplifying approximation, we model condenser and reboiler as equilibrium stages, consistent with the literature \citep{Akula.2021, Freguia.2003}. We assume Henry's law for \ce{CO2} \citep{Hartono.2014} and Raoult's law for water \citep{Perry.2008}. The cross heat exchanger is modeled as a plate-and-frame heat exchanger using the modeling methodology of \citet{akula.2019}, whereas the reboiler, condenser, cooler, and intercooler are sized from their heat duties assuming constant overall heat-transfer coefficients. 

The unit operation models include species and energy balances. In addition, we model the \ce{MEA}-\ce{CO2} reaction chemistry using the equilibrium reactions reported by \citet{Freguia.2002}, and impose these reactions in the condenser, reboiler, and cross heat exchanger. The absorber and desorber column models include the equilibrium chemistry extended by an approach-to-equilibrium model. Further details on the unit-operation models and property models are provided in \Crefrange{SI:unitoperations}{SI:properties}.
     
\subsection{Approach-to-equilibrium model}
The absorber and desorber are modeled as axially discretized columns, with material and energy balances and reaction chemistry solved in each segment. However, unlike the other unit operations, vapor–liquid equilibrium cannot be assumed in each column segment because absorption and desorption are limited by mass transfer and reaction kinetics. Rate-based models are the established approach for representing these rate limitations \citep{Zhang.2009}. However, rate-based models are prone to numerical instability when embedded in equation-oriented optimization. Rather than using a rate-based model, we represent rate-limitations through a calibrated approach-to-equilibrium function $\eta(\boldsymbol{\zeta}_{s,j};\boldsymbol{\beta})$. This function determines how far the vapor composition in each segment $j$ progresses toward its equilibrium value \cref{eq:approach_to_equilibrium}, where $\eta=1$ indicates equilibrium. This modeling approach reduces problem complexity in comparison to a more rigorous rate-based model to obtain well-behaved model for optimization:
\begin{equation}
    y_{\ce{CO2},s,j+1} = y_{\ce{CO2},s,j} - \eta(\boldsymbol{\zeta}_{s,j};\boldsymbol{\beta})\, \big( y_{\ce{CO2},s,j}- y_{\ce{CO2},s,j}^* \big),
    \qquad \forall\, s\in\mathcal{S},\, j\in\mathcal{J}.
    \label{eq:approach_to_equilibrium}
\end{equation}
Here, $y_{\ce{CO2},s,j}$ is the \ce{CO2} mole fraction in the vapor entering segment $j$ of the column in scenario $s$, while $y^{*}_{\ce{CO2},s,j}$ is the mole fraction at equilibrium, and $\mathcal{J}$ is the set of column segments.
The vector $\boldsymbol{\zeta}_{s,j}$ includes the local conditions of segment $j$: 
the vapor-to-liquid molar flow ratio $V_{s,j}/L_{s,j}$,
inlet temperature $T_{s,j}$,
pressure $p_{s,j}$,
solvent \ce{CO2} loading $\alpha_{s,j}$,
\ce{MEA} mass fraction in \ce{CO2}-free inlet solvent, $x^{\ce{MEA},\mathrm{lean}}_{s,j}$, and
superficial gas velocity $u^{\mathrm{v}}_{s,j}$. The vector $\boldsymbol{\beta}$ includes the coefficients of the model, which are calibrated from data as detailed below. The form of the function $\eta(\cdot)$ is chosen as
\begin{equation}
    \eta(\boldsymbol{\zeta}_{s,j};\boldsymbol{\beta}) = 1 - e^{-k(\boldsymbol{\zeta}_{s,j};\boldsymbol{\beta})\, \Delta z}, \qquad \forall\, s\in\mathcal{S},\, j\in\mathcal{J}.
    \label{eq:calibrationModel}
\end{equation}
Here, $\Delta z$ is the segment height resulting from the axial discretization. The function $k(\cdot\,;\boldsymbol{\beta})$ lumps the effects of gas- and liquid-phase mass transfer resistance, chemical reaction rates, interfacial area and hydrodynamics, and modeling approximations, and has the form
\begin{equation}
    \ln k(\boldsymbol{\zeta};\boldsymbol{\beta}) = \beta_0 + \sum_{n\in\mathcal{N}}{\beta_{n}} \, \zeta_{n} + \sum_{(n,m)\in\mathcal{M}}{\beta_{n,m}\, \zeta_n\,\zeta_m},
    \label{eq:fittingFunction}
\end{equation}
where $\mathcal{N}$ is the index set of calibration variables and $\mathcal{M} \subseteq\mathcal{N}\times\mathcal{N}$ is the set of index pairs of variables included to capture interactions, and $\boldsymbol{\beta}$ is the vector of calibration coefficients $\beta_0$, $\beta_n$, and $\beta_{n,m}$. 

The coefficients $\beta_0$, $\beta_n$, and $\beta_{n,m}$ are calibrated by fitting \cref{eq:fittingFunction} to data from simulations or experiments covering the relevant operating ranges. The coefficients are fitted separately for the absorber and the desorber columns. The coefficient values, ranges of data used for the fit, and prediction errors are reported in \Cref{SI:approachtoequilibrium}.

\subsection{Capture-rate targets}
A capture rate target is included as a constraint in the design problem. We specify a minimum capture rate either for every scenario,
\begin{equation}
    \dot{V}^{\mathrm{feed}}_{s} \, y^{\mathrm{feed}}_{\ce{CO2},s}
    -
    \dot{V}^{\mathrm{abs},\mathrm{out}}_{s} \, y_{\ce{CO2},s}^\mathrm{out}
    \ge
    \gamma^{\mathrm{cap}}\,
    \dot{V}^{\mathrm{feed}}_{s} \, y^{\mathrm{feed}}_{\ce{CO2},s},
    \qquad \forall\, s \in \mathcal{S},
    \label{eq:capturerate}
\end{equation}
or as a scenario-weighted average,
\begin{equation}
    \sum_{s\in\mathcal{S}} w_s\left(\dot{V}^{\mathrm{feed}}_{s} \, y^{\mathrm{feed}}_{\ce{CO2},s}
    -
    \dot{V}^{\mathrm{abs},\mathrm{out}}_{s} \, y_{\ce{CO2},s}^\mathrm{out} \right)
    \ge
    \gamma^{\mathrm{cap}}\,
    \sum_{s\in\mathcal{S}} w_s \left(\dot V^{\mathrm{feed}}_{s} \, y^{\mathrm{feed}}_{\ce{CO2},s}\right).
    \label{eq:capturerate_average}
\end{equation}
Here, $\dot{V}^{\mathrm{feed}}_{s}$ and $\dot{V}^{\mathrm{abs},\mathrm{out}}_{s}$ are the molar flows of the flue gas feed and the absorber outlet, respectively, with the \ce{CO2} mole fractions $y_{\ce{CO2},s}^\mathrm{feed}$ and $y_{\ce{CO2},s}^\mathrm{out}$, respectively. The parameter $\gamma^{\mathrm{cap}}$ is the capture rate target.

\subsection{Cost model}
The total annualized cost is selected as the objective function \cref{eq:objective} and combines annualized capital costs and scenario-weighted operating costs. The capital costs are estimated based on techno-economic models presented in \citet{Seider.2017} and the capital cost estimation methodology of the National Energy Technology Laboratory \citet{Theis.2021}. We assume that the equipment is made of stainless steel and escalated to 2023 using the Chemical Engineering Plant Cost Index. 

The bare erected costs of the columns include the column shell and its internals and are estimated based on column height and diameter. The bare erected cost of the cross heat exchanger is estimated using cost correlations for plate-and-frame heat exchangers in alignment with our process model, while the costs of reboiler, condenser, cooler, and intercooler are estimated based on shell-and-tube heat exchangers. We assume variable-speed pumps for the costing of the lean and rich solvent pumps. For the pumps, we include the energy demand for overcoming hydrostatic pressure and pressure changes between the absorption and desorption phases, assuming a fixed pump efficiency of \SI{65}{\percent}. The bare erected cost of the compressor is estimated assuming a multistage centrifugal compressor. The compressor power and cooling loads are estimated based on \citet{Kim.2025}. 

The bare erected costs of all components are converted to total overnight costs using cost factors following the methodology of \citet{Cvetic.2023,Theis.2021}. The total overnight cost in turn is annualized using the capital charge factor reported in \citet{Theis.2021}.

The operating costs include the utility costs for steam for solvent regeneration, electricity for compressor and pump operation, and cooling water with prices from \citet{Seider.2017}. We further include fixed operating costs for labor, maintenance, taxes and insurance \citep{Cvetic.2023}. The cost of \ce{MEA} makeup to compensate for degradation is estimated with prices from \citet{Haslbeck.2010}. Finally, other variable operating costs, including waste treatment and chemicals are based on \citet{Cvetic.2023} and scaled with the amount of \ce{CO2} captured. Further details on the cost models are included in \Cref{SI:cost}.

\section{Case study}
\label{sec:casestudy}
The goal of our case study is to characterize the cost benefits of accounting for part-load operation in the design of carbon capture processes. To this end, we solve the design and operation problem \cref{eq:problem} of the carbon capture process (\Cref{fig:processflowdiagram}) for flue gases from coal and \gls{ngcc} power plants, and compare designs obtained with and without accounting for part-load conditions. For simplicity, we assume an end-of-pipe carbon capture process that does not impact the operation or performance of the power plant. 

\subsection{Design cases}
We quantify the value of applying stochastic programming for design by comparing four cases: We consider a nominal case where we design the carbon capture plant, assuming a single scenario at the design operating point with \SI{80}{\percent} availability (\textit{nominal}). This reflects common design practice in industry, where availabilities of \SIrange{80}{85}{\percent} are commonly assumed \citep{Cvetic.2023, Schmitt.2022, Fisher.2007}, even though lower capacity factors lead to an underestimation of capture costs \citep{ElectricPowerResearchInstitute.2021}. The operating costs are evaluated by optimizing the operation for the fixed design. Next, we consider a multi-load design where we consider the design operating point along with multiple part-load scenarios and respective scenario weights (\textit{multi-load}). In both the \textit{nominal} and the \textit{multi-load} cases, we assume a \SI{95}{\percent} capture rate that must be achieved in all operating scenarios \cref{eq:capturerate}. Further, we include two additional cases (\textit{nominal (avg.)} and \textit{multi-load (avg.)}), which differ from the first two cases only by the assumption that the \SI{95}{\percent} capture rate needs to be achieved on average across all operating scenarios \cref{eq:capturerate_average}, introducing a scenario coupling constraint.

\subsection{Scenario definition}
\begin{figure}[!t]
    \centering
    \includegraphics[width=0.48\textwidth]{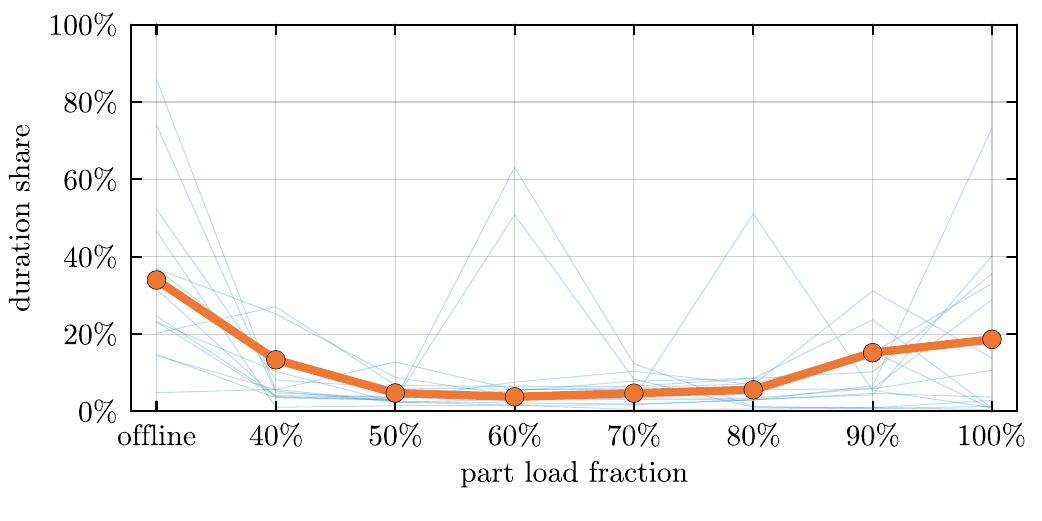}
    \caption{Duration share in part-load operation of sub-bituminous coal power plants in Texas in 2025. The medoid of the profiles is selected as the representative operation profile.}
    \label{fig:coal_partload}
\end{figure}

The scenarios $s\in\mathcal{S}$ differ in the flue gas inlet composition and mass flow. 
In the case study considering flue gases from a coal power plant, we combine data from various sources to construct the scenarios: For the coal case study, our scenarios represent operation at loads between \SI{100}{\percent} and \SI{40}{\percent} in \SI{10}{\percent} increments, resulting in seven operating scenarios. The scenario weights represent the duration of operation at each load point based on representative operating data of a coal power plant. For deriving the duration at each operational load level, we use operating load data from a fleet of sub-bituminous coal power plant units in Texas in 2025 \citep{U.S.EnvironmentalProtectionAgency.2026}. We allocate each operating point of the operating load data to the nearest \SI{10}{\percent}-part-load bin, to determine the duration share in each part-load operating point. We then select the medoid of all coal power plant units analyzed as the representative operation profile for the scenario weights of our case study (\Cref{fig:coal_partload}). 

The operating load data does not include flue gas composition and how it varies with load.  Therefore, we use flue gas flows and \ce{CO2} concentration at full load and part-load from simulations by \citet{Roeder.2014}, representative of a power plant of approximately \SI{600}{\mega\watt} capacity. We completed the flue gas composition data by assuming constant water and oxygen concentration in the pulverized coal power plant flue gas based on \citet{Schmitt.2022}. 
For the selected scenario weights, flue gas mass flow, and \ce{CO2} concentration, a total of \SI{1.8}{\mega\tonne\per\year} of \ce{CO2} are captured at a capture rate of \SI{95}{\percent}.

\begin{figure}[!t]
    \centering
    \begin{subfigure}[t]{0.48\textwidth}
        \centering
        \includegraphics[width=\linewidth]{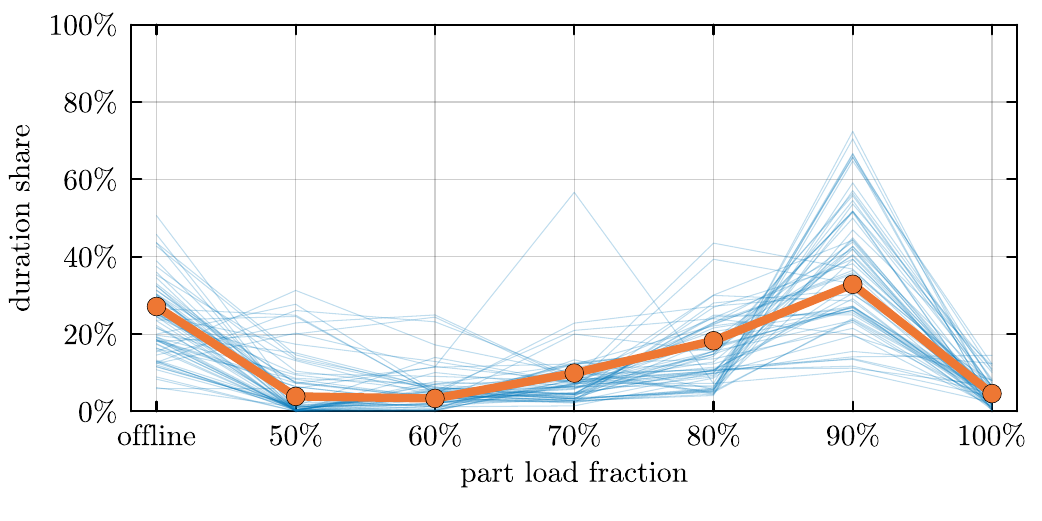}
        \caption{Cluster of plants with high capacity factor.}
        \label{fig:NGCC_high}
    \end{subfigure}
    
    \vspace{0.5em}
    
    \begin{subfigure}[t]{0.48\textwidth}
        \centering
        \includegraphics[width=\linewidth]{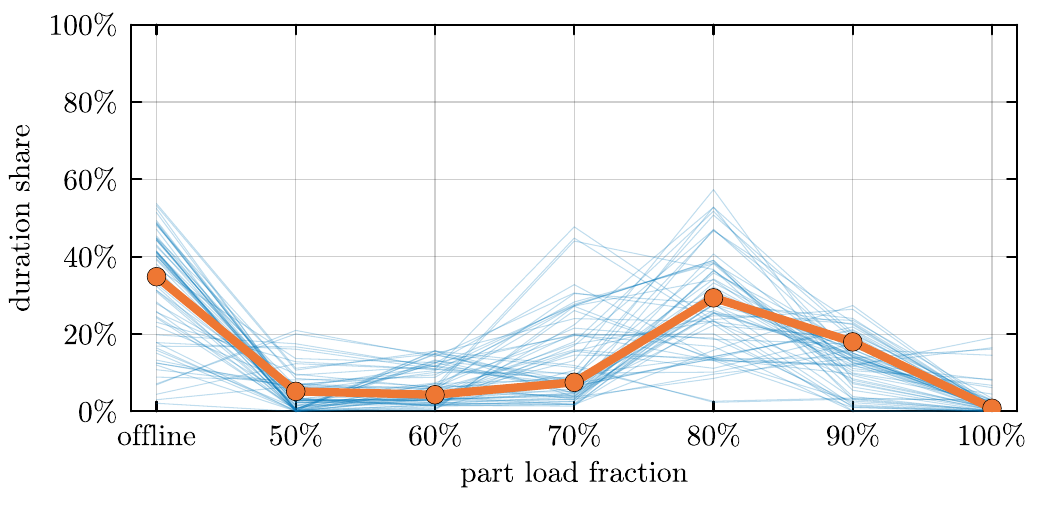}
        \caption{Cluster of plants with medium capacity factor.}
        \label{fig:NGCC_medium}
    \end{subfigure}
    
    \vspace{0.5em}
    
    \begin{subfigure}[t]{0.48\textwidth}
        \centering
        \includegraphics[width=\linewidth]{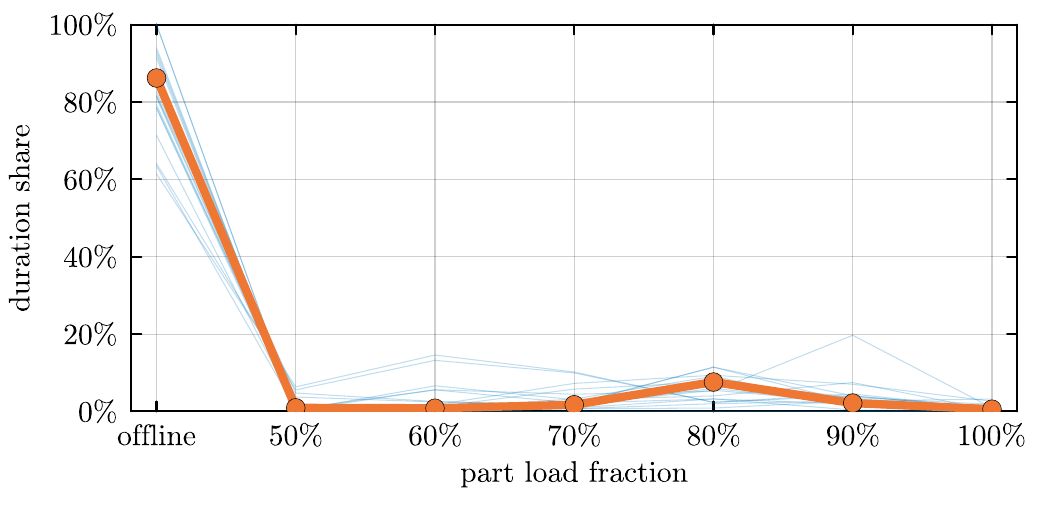}
        \caption{Cluster of plants with low capacity factor.}
        \label{fig:NGCC_low}
    \end{subfigure}
    \caption{Duration share in part-load operation of natural gas combined cycle power plants in Texas in 2025. Representative operation profiles determined via k-medoids clustering are highlighted in orange.}
    \label{fig:NGCC_partload}
\end{figure}

For the \gls{ngcc} case study, our scenarios represent operation at loads between \SI{100}{\percent} and \SI{50}{\percent} in \SI{10}{\percent} increments, resulting in six operating scenarios. Again, the duration of operation at each operating point is based on representative operating data of an \gls{ngcc} unit from Texas in 2025 \citep{U.S.EnvironmentalProtectionAgency.2026}. We again determine the duration share in each part-load operating point for each \gls{ngcc} unit. Due to the relatively large number of 163 units, we apply k-medoids clustering \citep{Kaufman.1987} to determine three clusters of \gls{ngcc} units, with high, medium, and low capacity factor (\Cref{fig:NGCC_partload}). We select the medoid of the high capacity factor cluster for the scenario weights in our case study.

Similar to the coal case study, we use simulated flue gas mass flows and \ce{CO2} concentrations at full load and part-load from \citet{Alcaraz-Calderon.2019}, for a power plant with a capacity of approximately \SI{700}{\mega\watt}. For the selected scenario weights, flue gas mass flow, and \ce{CO2} concentration, a total of \SI{1.4}{\mega\tonne\per\year} of \ce{CO2} are captured at a capture rate of \SI{95}{\percent}.

Due to the large flue gas mass flows reported for both coal and \gls{ngcc} plants, we assume multiple identical carbon capture plants to be installed in parallel to process the flue gases. Specifically, we assume two carbon capture trains in the coal case study and four trains in the \gls{ngcc} case study, corresponding to a flue gas mass flow of approximately \SI{1000}{\tonne\per\hour} per capture train at full load in both case studies. Installing multiple capture trains is an option to limit equipment size and increase both redundancy and operational flexibility \citep{Agbonghae.2014}. For simplicity, we assume the flue gas stream to be evenly distributed between capture trains, neglecting the option of selective shutdowns of trains and of shared equipment between trains.

\subsection{Approach-to-equilibrium model calibration}
We calibrate the approach-to-equilibrium model \cref{eq:fittingFunction} for the absorption and desorption columns by fitting the calibration coefficients $\beta_0$, $\beta_n$, and $\beta_{n,m}$ to data simulated using the Aspen Plus \textit{CCSI steady state MEA model} \citep{U.S.DepartmentofEnergy.2022} for a wide range of operating conditions (\Cref{tab:operatingrange}). We tailor the range of operating conditions to each case study by adjusting the inlet flue gas \ce{CO2}-concentration to a representative range.
\begin{remark}
This calibration emulates the industrial practice of calibrating equilibrium-based models against operational data. Since operational data are not available to us, we instead generate the calibration data with the rate-based Aspen Plus \textit{CCSI steady state MEA model} \citep{U.S.DepartmentofEnergy.2022}.
\end{remark} 

For the coal case study, we generate 360,000 simulated data points sampled from the range of operating conditions. The approach-to-equilibrium model results in an average prediction error of \SI{1.5}{\percent} for the column segment output \ce{CO2} mole fraction, with the 5th to 95th percentile prediction error ranging from \SIrange{-3}{6.6}{\percent} (\Cref{fig:predictionerror_coal}).
For the \gls{ngcc} case study, we generate 190,000 simulated data points, obtaining an average prediction error of \SI{0.7}{\percent}. The 5th to 95th percentile prediction error is wider than in the coal power plant case study, ranging from \SIrange{-6.3}{10.8}{\percent} (\Cref{fig:predictionerror_gas}).

We validate the carbon capture process modeling approach by simulating demonstrator carbon capture processes from the literature under various operating points and comparing simulated to observed capture rates.
In each simulated operating point, we fix inlet flue gas and lean solvent compositions and simulate the carbon capture rate, with the plant design fixed to the specifications of the demonstrator plants as given in the literature. 

For validating the model for the coal case study, we simulate the operation of a carbon capture plant using \SI{30}{\wtpercent} aqueous MEA based on the operating conditions during a measurement campaign considering flue gases with a \ce{CO2} concentration of \SIrange{9}{12}{\molpercent}, based on \citet{Morgan.2018, Akula.2021}. For validating the model for the \gls{ngcc} case study, we in turn simulate the conditions of a measurement campaign for treating flue gases with a lower \ce{CO2} concentration of \SIrange{3.4}{3.7}{\molpercent}, based on \citet{Oi.2021, Gjernes.2017, Faramarzi.2017, Hamborg.2014}.

\begin{figure}[!t]
\centering
\begin{subfigure}[t]{0.48\textwidth}
    \centering
    \includegraphics[width=\linewidth]{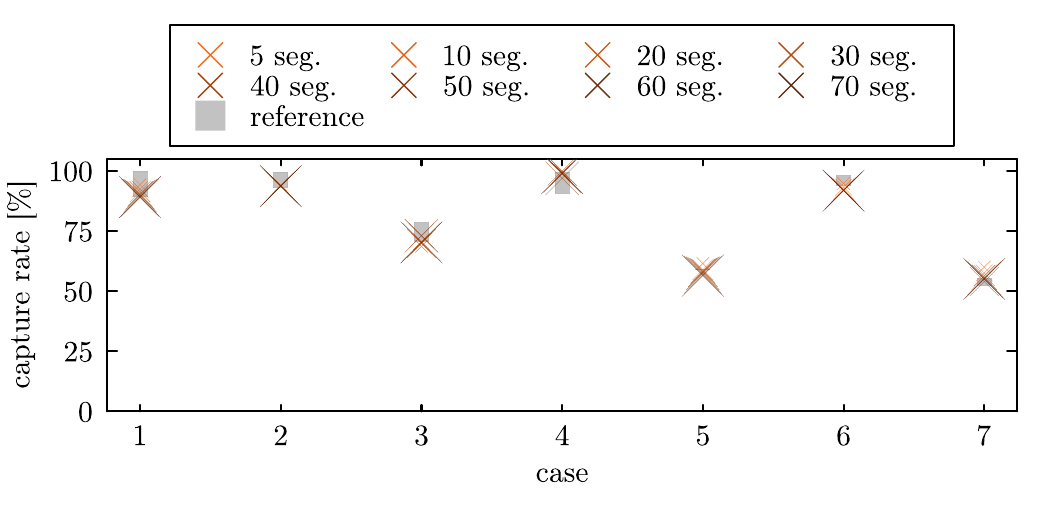}
    \caption{Validation results for flue gas compositions representative of coal-fired power plants. The inlet flue gas and lean solvent conditions of the measurement campaign are reported by \citet{Morgan.2018} and summarized by \citet{Akula.2021}.}
    \label{fig:validation_coal}
\end{subfigure}

\vspace{1.5em}

\begin{subfigure}[t]{0.48\textwidth}
    \centering
    \includegraphics[width=\linewidth]{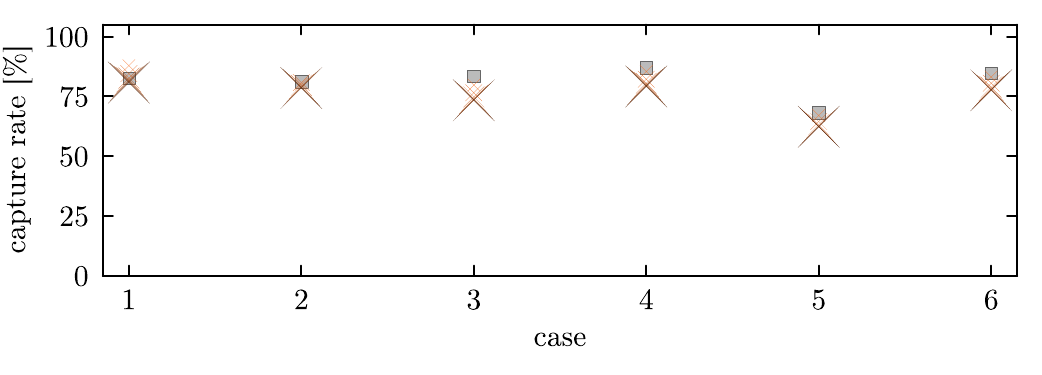}
    \caption{Validation results for flue gas compositions representative of natural gas combined cycle plants. The inlet flue gas and lean solvent conditions of the measurement campaign are provided in \Cref{tab:validation_gas_fluegas,tab:validation_gas_solvent}.}
    \label{fig:validation_gas}
\end{subfigure}

\vspace{1em}
\caption{Simulated capture rates ($\times$) in comparison to measurement results from demonstrator plants (reference). Simulation results are obtained by fixing the compositions and flow rates of inlet flue gas and solvent to the various conditions during measurements of the demonstrator plants. The simulation results are shown for a a number of segments (seg.) ranging from 5 to 70.}
\label{fig:validation}
\end{figure}

The validation demonstrates agreement between the simulated and measured capture rate, with the number of segments specified in our model having only a minor impact on the results (\Cref{fig:validation}).

\subsection{Solution strategy}
Solving the nonlinear design problem requires a careful initialization strategy. We obtain good starting values for the nonlinear optimization problem by first simulating the carbon capture process without recycles. To this end, we tear the lean solvent recycle stream, and fix the lean solvent compositions and mass flows, along with the plant design and reboiler duty using reasonable values. We simulate the model without recycles for a single scenario and select a coarse resolution for the columns of 10 segments. The simulation results are used as starting values for the design optimization problem with the recycle stream restored and the same column resolution and number of scenarios. The column resolution and the number of scenarios are then successively increased, initialized using the previous results and applying linear interpolation between column segments if necessary. Our design results are based on a resolution of 30 segments for both absorber and desorber columns, which is sufficiently high for the convergence of the capture rate (\Cref{fig:validation}). After obtaining a feasible solution for all part-load scenarios at the target column resolution, we re-solve the design optimization 20 times using this solution for initialization. In each run, the initial design variables are randomly perturbed within a \SI{20}{\percent} range to explore nearby local optima. We report the design with the lowest total annualized cost. However, in practice, the perturbations yield only marginal cost reductions.

The design optimization problem is modeled in \textit{Julia} with \textit{JuMP} v1.29.3 \citep{Lubin.2023} and solved with \textit{IPOPT} v1.13.0 \citep{Wachter.2006} using the linear solver \textit{MA57} \citep{Duff.2004}.

\section{Results}
\label{sec:results}
We first compare the designs obtained for the coal and \gls{ngcc} case studies, before analyzing the resulting costs and their sensitivities. All results reported in this section can be reproduced using the scripts provided in our repository.\footnote{\url{https://github.com/mit-shin-group/part-load-design-CCS-2026}}

\begin{figure}[!t]
    \centering
    \begin{subfigure}[t]{\columnwidth}
        \centering
        \includegraphics[width=0.48\columnwidth,trim={0 0.5cm 0 0},clip]{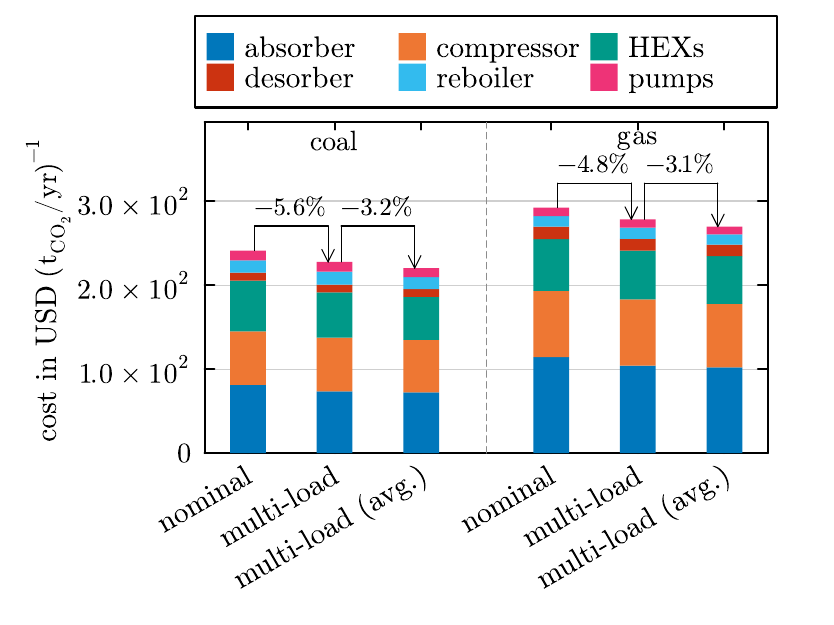}
        \caption{Total plant cost and main contributors. The designs of the two nominal cases are identical.}
        \label{fig:capex}
    \end{subfigure}

    \vspace{0.5em}
    
    \begin{subfigure}[t]{\columnwidth}
        \centering
        \includegraphics[width=0.48\columnwidth,trim={0 0.5cm 0 0},clip]{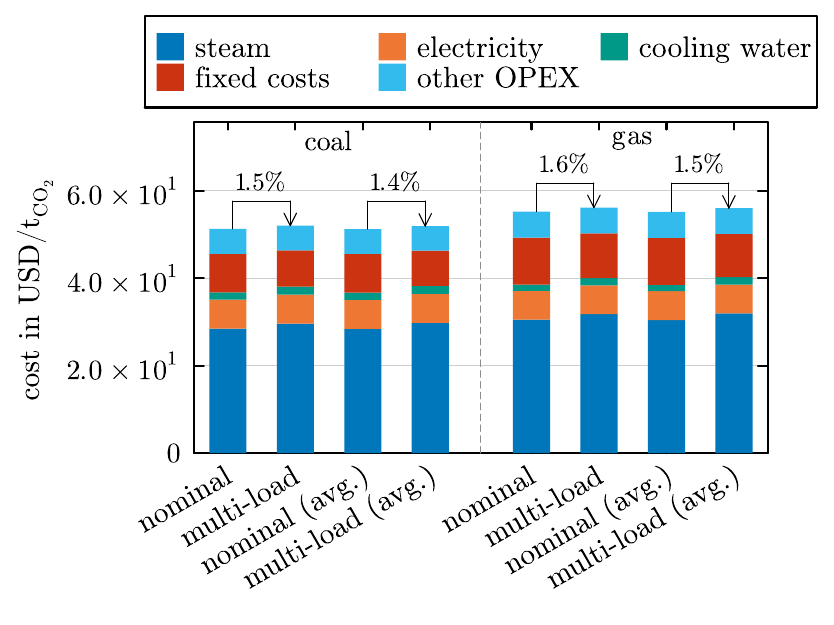}
        \caption{Annual operating cost of the carbon capture process and main contributors.}
        \label{fig:opex}
    \end{subfigure}
    
    \vspace{0.5em}
    
    \begin{subfigure}[t]{\columnwidth}
        \centering
        \includegraphics[width=0.48\columnwidth,trim={0 0.5cm 0 0},clip]{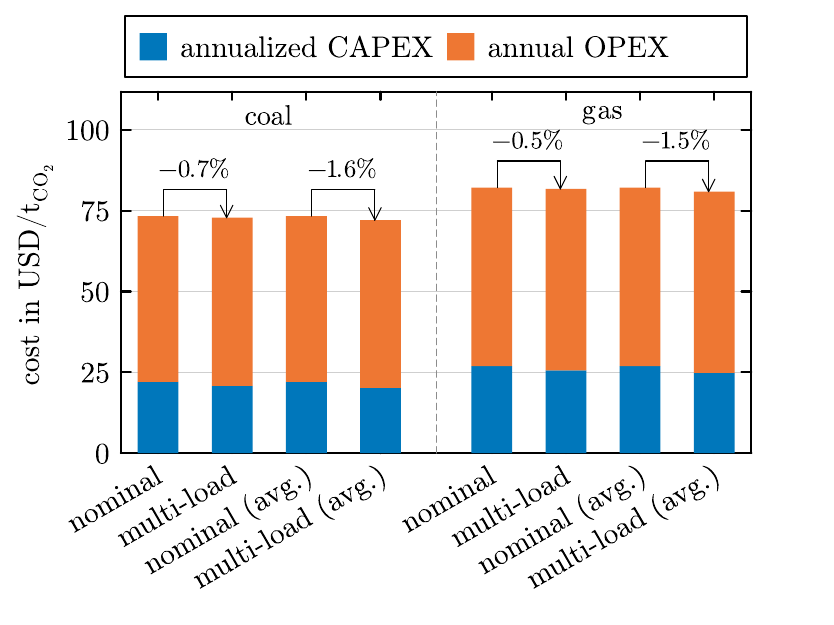}
        \caption{Total annualized cost of carbon capture.}
        \label{fig:TAC}
    \end{subfigure}
    \caption{Costs of the capture process for nominal and multi-load designs in the coal (left) and natural gas (right) case studies, at a \SI{95}{\percent} capture rate per load scenario or on average (avg.) across scenarios. \newline\textit{HEX: heat exchanger, OPEX: operating expenditures}}
    \label{fig:cost}
\end{figure}

Even though both case studies consider the same flue gas mass flows`per capture train, most units of the carbon capture system designed for the coal power plant are larger than in the \gls{ngcc} case study (\Cref{tab:componentsizes}). This is due to the higher \ce{CO2} mole fraction in the flue gas of the coal power plant, which results in larger solvent and \ce{CO2}-product mass flows. Despite the larger equipment sizes, the coal capture plant achieves \SIrange{17}{18}{\percent} lower total plant costs than in the \gls{ngcc} case study when normalized by capture capacity, due to economies of scale (\Cref{fig:capex}). Similarly, the higher \ce{CO2} mole fraction results in \SI{7}{\percent} lower specific operating costs for the coal capture process compared to the capture process of the \gls{ngcc} plant in all cases (\Cref{fig:opex}). 
\begin{table*}[!t]
\centering
\caption{Process equipment sizes optimized assuming a single nominal operating point (nominal) and relative change in designs accounting for multiple part-load scenarios (multi-load) and a scenario-averaged capture rate (avg.).\\ \textit{abbreviations: absorber (abs.), desorber (des.), heat exchanger (HEX)}}
    \begin{tabular}{lrrr}
    \toprule
    variable & \begin{tabular}[c]{@{}r@{}}nominal\end{tabular} & multi-load & \begin{tabular}[c]{@{}r@{}}multi-load\\ (avg.)\end{tabular} \\
    \midrule
    \multicolumn{4}{l}{\textit{coal case study}} \\
    \midrule
    packed height (abs.) & $24.3\,\si{\meter}$              & $-10.5\,\si{\percent}$    & $-9.9\,\si{\percent}$ \\
    diameter (abs.) & $11.9\,\si{\meter}$                   & $+1.1\,\si{\percent}$     & $-0.5\,\si{\percent}$ \\
    packed height (des.) & $5.0\,\si{\meter}$               & $-11.3\,\si{\percent}$    & $-8.5\,\si{\percent}$ \\
    diameter (des.) & $7.9\,\si{\meter}$ &                  $+2.2\,\si{\percent}$       & $-1.1\,\si{\percent}$ \\
    condenser area & $6.6\times10^{3}\,\si{\square\meter}$  & $-4.0\,\si{\percent}$     & $-9.3\,\si{\percent}$ \\
    cooler area & $3.7\times10^{3}\,\si{\square\meter}$     & $+17.1\,\si{\percent}$    & $+4.2\,\si{\percent}$ \\
    intercooler area & $4.4\times10^{3}\,\si{\square\meter}$ & $-18.0\,\si{\percent}$   & $-24.4\,\si{\percent}$ \\
    reboiler area & $3.9\times10^{3}\,\si{\square\meter}$   & $+4.5\,\si{\percent}$     & $-1.1\,\si{\percent}$ \\
    compressor power & $17.8\,\si{\mega\watt}$              & $0.0\,\si{\percent}$      & $-2.9\,\si{\percent}$ \\
    lean pump flow & $1.9\times10^{4}\,$\si{\gpm}                 & $+4.5\, \si{\percent}$    & $-3.6\,\si{\percent}$ \\
    rich pump flow & $2.0\times10^{4}\,$\si{\gpm}                 & $+4.2\,\si{\percent}$     & $-3.5\,\si{\percent}$ \\
    \midrule
    \multicolumn{4}{l}{\textit{natural gas case study}} \\
    \midrule
    packed height (abs.) & $15.8\,\si{\meter}$              & $-9.6\,\si{\percent}$     & $-10.3\,\si{\percent}$ \\
    diameter (abs.) & $10.8\,\si{\meter}$                   & $+0.2\,\si{\percent}$     & $-0.5\,\si{\percent}$ \\
    packed height (des.) & $6.7\,\si{\meter}$               & $-11.0\,\si{\percent}$     & $-7.7\,\si{\percent}$ \\
    diameter (des.) & $4.7\,\si{\meter}$                    & $+1.9\,\si{\percent}$     & $-2.2\,\si{\percent}$ \\
    cross HEX area & $1.4\times10^{4}\,\si{\square\meter}$  & $-28.5\,\si{\percent}$    & $-24.3\,\si{\percent}$ \\
    condenser area & $2.1\times10^{3}\,\si{\square\meter}$  & $-0.6\,\si{\percent}$     & $-7.4\,\si{\percent}$ \\
    cooler area & $7.9\times10^{2}\,\si{\square\meter}$     & $+19.1\,\si{\percent}$    & $+4.7\,\si{\percent}$ \\
    intercooler area & $1.5\times10^{3}\,\si{\square\meter}$ & $-1.9\,\si{\percent}$    & $-2.5\,\si{\percent}$ \\
    reboiler area & $1.3\times10^{3}\,\si{\square\meter}$   & $+4.0\,\si{\percent}$     & $-4.6\,\si{\percent}$ \\
    compressor power & $5.7\,\si{\mega\watt}$               & $0.0\,\si{\percent}$      & $-7.1\,\si{\percent}$ \\
    lean pump flow & $7.0\times10^{3}\,$\si{\gpm}                 & $+3.4\,\si{\percent}$     & $-4.1\,\si{\percent}$ \\
    rich pump flow & $7.4\times10^{3}\,$\si{\gpm}                 & $+3.2\,\si{\percent}$     & $-4.3\,\si{\percent}$ \\
    \bottomrule
    \end{tabular}
    \label{tab:componentsizes}
\end{table*}
When comparing designs within each case study, the units are predominantly smaller in the process designs resulting from the two multi-load cases compared to the nominal case  (\Cref{tab:componentsizes}). 
In the coal power case study, absorber and desorber packing heights are reduced by \SIrange{9}{11}{\percent}, and cross heat exchanger surface area is reduced by up to \SI{50}{\percent} compared to the nominal reference. Similarly, in the \gls{ngcc} case study, packing heights are reduced by \SIrange{8}{11}{\percent}, and cross heat exchanger surface area is reduced by up to \SI{29}{\percent}.

Smaller equipment sizes in the multi-load cases translate into lower total plant costs, with costs being lowest for the multi-load case with a scenario-weighted average capture target. Total plant costs in the coal case study are reduced substantially by \SIrange{6}{9}{\percent} and in the \gls{ngcc} case study by \SIrange{5}{8}{\percent} (\Cref{fig:capex}). These reductions are driven by smaller absorber columns and heat exchangers, which contribute more than half of the capital cost in the nominal cases.

The total plant cost reductions in the multi-load cases are partially offset by higher annual operating costs. The operating costs are dominated by the regeneration steam demand, which increases by \SIrange{4}{5}{\percent} in the multi-load cases (\Cref{fig:opex}) due to the smaller equipment sizes, in particular of the cross heat exchanger. However, the increase in operating cost is outweighed by the reduction of the capital cost of the cross heat exchanger and columns.

Overall, the multi-load cases result in a modest overall reduction in total annualized cost by \SIrange{0.7}{1.6}{\percent} for the coal power plant case study and by \SIrange{0.5}{1.5}{\percent} for the \gls{ngcc} case study (\Cref{fig:TAC}).

The reductions in total annualized cost are larger in the case considering a scenario-weighted average capture rate. In the multi-load case with scenario-weighted capture rate, \ce{CO2} capture is shifted from full-load to part-load periods: In the coal case study, a capture rate of only \SI{92}{\percent} is achieved at full load, while the capture rate reaches up to \SI{98}{\percent} at part load in the multi-load case (\Cref{fig:capture_rate_vs_part_load_coal}). 
The shift is even more pronounced for the \gls{ngcc} case study, where capture rates vary between \SIrange{88}{98}{\percent} from full load to part load (\Cref{fig:capture_rate_vs_part_load_gas}). The nominal cases also shift \ce{CO2} capture to part-load operation, although to a lesser degree, and the total annualized cost of the nominal case with scenario-weighted average capture rate is reduced only marginally.

\subsection*{Sensitivities}
\begin{figure*}[!ht]
    \centering
    \begin{subfigure}{0.48\textwidth}
        \centering
        \includegraphics[width=0.98\textwidth,trim={0 3cm 0 0},clip]{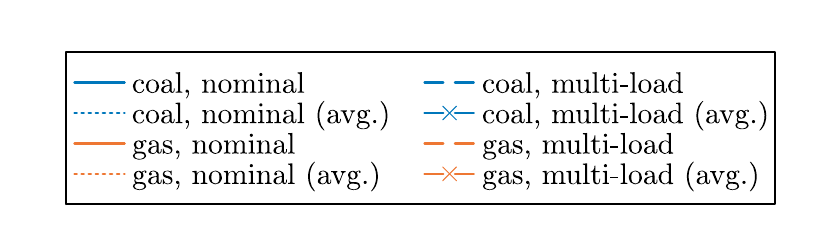}
    \end{subfigure}

    \begin{subfigure}[t]{0.48\textwidth}
        \centering
        \includegraphics[width=0.98\textwidth,trim={0 0.5cm 0 0},clip]{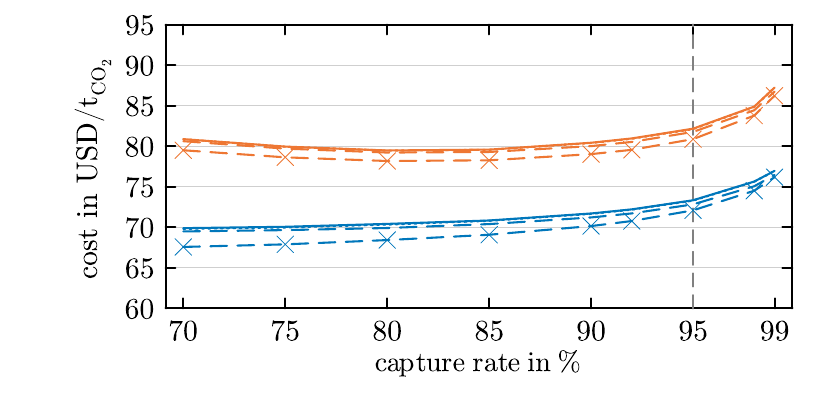}
        \caption{Sensitivity to variations in capture rate target.}
        \label{fig:sensitivity_capture_rate}
    \end{subfigure}
    \hfill
    \begin{subfigure}[t]{0.48\textwidth}
        \centering
        \includegraphics[width=0.98\textwidth,trim={0 0.5cm 0 0},clip]{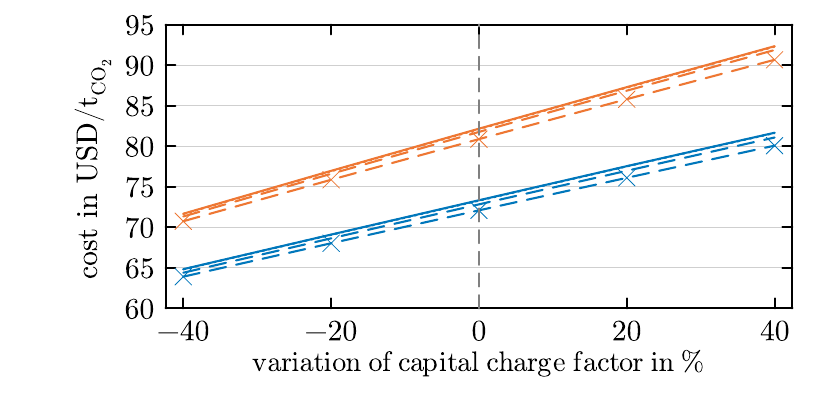}
        \caption{Sensitivity to variations of capital charge factor.}
        \label{fig:sensitivity_CCF}
    \end{subfigure}
    
    \vspace{0.3em}
    
    \begin{subfigure}[t]{0.48\textwidth}
        \centering
        \includegraphics[width=0.98\textwidth,trim={0 0.5cm 0 0},clip]{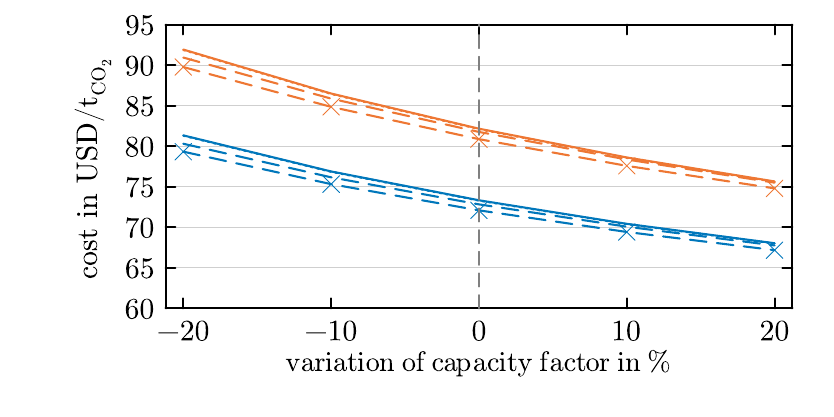}
        \caption{Sensitivity to variations in capacity factor.}
        \label{fig:sensitivity_flh}
    \end{subfigure}
    \hfill
    \begin{subfigure}[t]{0.48\textwidth}
        \centering
        \includegraphics[width=0.98\textwidth,trim={0 0.5cm 0 0},clip]{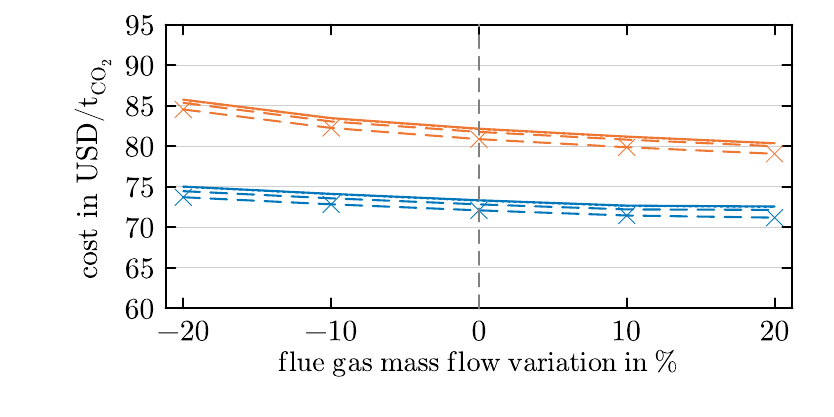}
        \caption{Sensitivity to variations in flue gas mass flow.}
        \label{fig:sensitivity_massflow}
    \end{subfigure}
    
    \caption{Sensitivity of the total annualized cost of carbon capture per tonne of \ce{CO2} catu. The baseline capacity factors are \SI{50.1}{\percent} and \SI{59.7}{\percent} for the coal and the natural gas combined cycle power plant, respectively (see \Cref{fig:NGCC_partload}). The capacity factor is varied by applying a scaling factor to all scenario weights, thereby reducing the amount of time the plant is offline. The baseline flue gas mass flows per capture train are \SI{282.3}{\kilo\gram\per\second} and \SI{285.5}{\kilo\gram\per\second} at full load for the coal and the natural gas combined cycle power plant, respectively (see \Cref{SI:part_load}). The baseline capital charge factor is \SI{7.7}{\percent} (see \Cref{SI:cost}).}
    \label{fig:sensitivity}
\end{figure*}

Higher shares of part-load operation reduce annual utilization and result in larger overall cost savings of multi-load designs, as operating conditions diverge further from the nominal design point. For example, total annualized cost savings increase from \SIrange{0.5}{1.5}{\percent} for \gls{ngcc} plants with high utilization (\Cref{fig:NGCC_high}) to \SIrange{1}{3}{\percent} (\Cref{fig:cost_gas_intermediate}) when considering a plant with intermediate utilization (\Cref{fig:NGCC_medium}). Similarly, decreasing the capacity factor by \SI{20}{\percent} increases the cost benefits of the multi-load cases to up to \SI{2.4}{\percent} and \SI{2.3}{\percent} in the \gls{ngcc} and the coal case studies, respectively (\Cref{fig:sensitivity_flh,fig:sensitivity_flh_costsaving}).

Certain technical and economic parameters affect the total annualized cost far more strongly than the switch from nominal to a multi-load design (\Cref{fig:sensitivity}): Reducing the target capture rate to \SI{80}{\percent} reduces the specific cost of carbon capture by \SIrange{3}{4}{\percent} in the nominal cases of the case studies (\Cref{fig:sensitivity_capture_rate}). Further reducing the capture rate results in an increase in the specific cost for the \gls{ngcc} cases, as the benefits of economies of scale are diminished with smaller equipment sizes. For the coal cases, the specific cost plateaus later at a capture target of \SI{70}{\percent}. Conversely, increasing the capture target from \SI{95}{\percent} to \SI{99}{\percent} raises the average specific cost of carbon capture by \SIrange{5}{6}{\percent}, indicating an expected accelerated cost increase with high capture rates.

The specific total cost of carbon capture increases almost linearly with the capital charge factor (\Cref{fig:sensitivity_CCF}), which would increase from the baseline value of \citet{Theis.2021} with a shorter lifetime of equipment or higher financing costs. An increase by \SI{40}{\percent}, corresponding to a return to levels assumed in previous reports \citep{Theis.2021}, would result in an increase in the specific cost by \SIrange{11}{12}{\percent} in the nominal cases.

In contrast, increases in both capacity factor and average flue gas mass flow reduce specific total annualized costs. In the nominal cases, a \SI{20}{\percent} increase in capacity factor reduces specific costs by \SIrange{7}{8}{\percent} (\Cref{fig:sensitivity_flh}), as the fixed capital cost is spread over a greater amount of \ce{CO2} captured. Increasing the flue gas mass flow by \SI{20}{\percent} yields a smaller reduction in specific total annualized cost of \SIrange{1}{2}{\percent} in the nominal cases (\Cref{fig:sensitivity_massflow}), driven by economies of scale. While these parameter variations shift the overall cost levels, the relative savings of the multi-load designs persist across all sensitivity cases (\Cref{fig:sensitivity_costsaving}). We therefore expect similar percentage savings under different operating conditions, cost parameters, and equipment sizes.

\section{Discussion}
\label{sec:discussions}
Accounting for part-load operation in the design stage of carbon capture plants has a large impact on equipment sizing: The multi-load cases feature predominantly smaller equipment than the nominal cases (\Cref{tab:componentsizes}) and substantially reduce capital costs (\Cref{fig:capex}). The resulting reductions in total annualized costs are moderate; however, they are achieved assuming no changes to solvent technology and are exclusively based on better sizing decisions in the design. Designing the carbon capture process based on a single design point and neglecting part-load operation misrepresents the tradeoff between capital and operational expenditures, resulting in higher overall costs in nominal designs.

The tradeoff is most visible in the cross heat exchanger, whose surface area is reduced by up to \SI{50}{\percent} (coal) and \SI{29}{\percent} (\gls{ngcc}) (\Cref{tab:componentsizes}). The smaller cross heat exchanger reduces the amount of internal heat regeneration and increases the need for both external heating and cooling. Thus, the operating costs increase driven by a \SIrange{4}{5}{\percent} higher regeneration steam demand (\Cref{fig:opex}). Accepting a higher specific energy demand in exchange for smaller cross heat exchanger surface area is only optimal once the limited duration of full-load operation is represented in the stochastic optimization problem. 
 
The lowest specific annualized costs across all case studies are achieved when accounting for multiple part-load scenarios while targeting a scenario-weighted average capture rate (\Cref{fig:TAC}, \textit{multi-load (avg.)}). Targeting an average capture rate provides additional operational flexibility that allows \ce{CO2} capture to be shifted from full-load to part-load periods, reducing peak mass flows in the capture process system. As a consequence, many equipment sizes are further reduced along with their capital costs. Only the cooler area remains larger than in the nominal design, as the substantial reduction in cross heat exchanger size reduces the amount of internal heat regeneration and increases the need for external heating and cooling supply. The increase in external energy demand is reflected in the component sizes of the reboiler and cooler in the multi-load case. However, in the case considering average capture rates (multi-load avg.), these increases in size are largely compensated (\Cref{tab:componentsizes}) by shifting capture to part-load periods. These resulting cost benefits are consistent across all cases of the sensitivity analysis (\Cref{fig:sensitivity}) and are achieved for the same annual \ce{CO2} emission reductions.

The additional flexibility can only be exploited when considering part-load operation at the design stage: Although the nominal cases also shift \ce{CO2} capture to part-load operation, no capital cost savings are achieved, as the nominal design is based on the design operating point. As the reduction in operating costs from shifting capture to part load is small compared to cost savings from capital costs, the total annualized costs are reduced only marginally in the nominal case with scenario-weighted average capture rate.

The cost savings achieved by applying the multi-load design approaches are modest compared to the cost changes resulting from plausible variations of key technical and economic parameters, such as the capture rate target, the capital charge factor, the capacity factor, and the flue gas mass flow (\Cref{fig:sensitivity}). For example, designing the plant with a capture rate target of \SI{80}{\percent} reduces both plant size and the specific capture cost by \SIrange{3}{4}{\percent} in the nominal cases (\Cref{fig:sensitivity_capture_rate}). However, the cost savings are persistent across all cases in the sensitivity analysis, indicating a consistent cost benefit of accounting for the part-load conditions in the plant design.

In this work, part-load operation is represented as a set of steady-state scenarios, neglecting ramping, start-up, and shutdown dynamics, which may affect the estimated operational flexibility of the capture plant. We further focus on the carbon capture plant and thus assume no changes to power plant operation or performance, neglecting interactions between the power and capture plant.
Future research should therefore include both the power plant and the carbon capture process in the design problem, to better represent the economic and environmental performance of the coupled system. The co-design problem can capture interactions, such as steam and electricity supply by the power plant and additionally enables process integration between the power plant and the capture process, offering further cost reduction potential. Accounting for dynamic operation of both plants would further enable assessing the flexibility that the coupled system can provide to the grid. However, the anticipated increase in model complexity needs to be addressed, e.g., via the application of GPU-accelerated optimization techniques \citep{Shin.2024}.

\section{Conclusions}
\label{sec:conclusions}
This study applied a two-stage stochastic optimization formulation to simultaneously optimize the design and operation of a solvent-based carbon capture process system considering multiple part-load operation profiles derived from real-world power plants. 

The process model in the stochastic optimization problem represents rate limitations in the absorption and desorption columns using an approach-to-equilibrium model to reduce model complexity in comparison with a more rigorous rate-based model. The process model further includes techno-economic models to estimate the total annualized cost of the carbon capture process while achieving capture targets.

Our results indicate that accounting for multiple part-load scenarios results in a noticeably different trade-off between capital and operating costs. As a consequence, the designs obtained following the multi-load approach feature predominantly smaller equipment, which translates into substantial reductions of the total plant costs by \SIrange{5}{9}{\percent} in the cases considered. The reductions in total plant cost are partially offset by increased operating costs, resulting in reductions of \SIrange{0.5}{1.6}{\percent} in total annualized cost compared to the nominal design. Although these savings may appear modest, the sensitivity analysis demonstrates that they persist in both direction and magnitude across variations of the key technical and economic parameters (\Cref{fig:sensitivity}). We therefore expect similar relative savings under different operating conditions, cost parameters, and equipment sizes.

This study demonstrates the cost-saving potential of designing carbon capture processes considering real-world part-load operating conditions via stochastic optimization. Part-load variability should therefore be considered at the design stage, particularly for plants that operate as load-following generators.

\section*{CRediT authorship contribution statement}
\textbf{David Yang Shu}: Conceptualization, Methodology, Software, Validation, Formal analysis, Investigation, Data Curation, Writing - Original Draft, Visualization, Project administration. \textbf{Boxun Huang}: Methodology, Software, Validation, Formal analysis, Investigation, Data Curation, Writing - Review \& Editing. \textbf{Yurim Kim}: Methodology, Validation, Investigation, Data Curation, Writing - Review \& Editing. \textbf{Rahul Gandhi}: Methodology, Validation, Writing - Review \& Editing. \textbf{Randall Field}: Methodology, Validation, Writing - Review \& Editing, Supervision. \textbf{Sungho Shin}: Conceptualization, Methodology, Software, Validation, Formal analysis, Resources, Writing - Review \& Editing, Supervision, Project administration, Funding acquisition.

\section*{Declaration of interests}
DS reports financial support was provided by the U.S. Department of Energy. RG reports a relationship with Shell Technology Center that includes: employment. 

\section*{Declaration of generative AI use}
During the preparation of this work the authors used Claude Sonnet 4.6 in order to debug and generate code, and Claude Fable 5 to check and polish the writing. After using these tools, the authors reviewed and edited the content as needed and take full responsibility for the content of the published article.

\section*{Acknowledgments}
DS, BH, and SS received support from the Future Energy Systems Center through the MIT Energy Initiative. YK received support from the Sejong Science Fellowship. This publication is based on work supported by the U.S. Department of Energy under contract DE-AC02-06CH11357. Funding was provided by the U.S. DOE Hydrocarbons and Geothermal Energy Office under FWP-DOE-145-SE-14205HPC4E-25. 

\sloppy
\printbibliography

\newpage
\appendix
\crefalias{section}{appendix}

\section{Modeling details} \label{SI:models}
This section provides the full model description of the carbon capture process system (\Cref{fig:processflowdiagram}). The transport and geological storage of \ce{CO2} are excluded from both the cost and process model.

\subsection{Equilibrium reactions}
We model \SI{30}{\wtpercent} aqueous \gls{mea}, neglecting the volatility of \gls{mea} and the solubility of \ce{N2} and \ce{O2} as a simplifying assumption. The liquid phase is described by the eleven species in \Cref{tab:species}. Speciation is governed by five equilibrium reactions:
\begin{reactions}
MEAH+ &<=> MEA + H+ \label{eq:reaction1}\\
HCO3- &<=> CO3^2- + H+ \\
H2O &<=> H+ + OH- \\
MEACOO- + H2O &<=> MEA + HCO3- \\
CO2 + H2O &<=> HCO3- + H+\label{eq:reaction5}
\end{reactions}
Equilibrium constants are expressed as $\ln K_r = \mathrm{A}_r + \mathrm{B}_r/T + \mathrm{C}_r\ln T + \mathrm{D}_r T$, with coefficients by \citet{Freguia.2002}. These reactions are imposed at equilibrium in every unit (absorber, desorber, cross heat exchanger, reboiler, condenser, cooler, and mixer).

\begin{table}[!ht]\centering
\caption{Species considered in the carbon capture process.}
\label{tab:species}
\begin{tabular}{p{0.15\columnwidth}p{0.75\columnwidth}}
\toprule
 & species \\
\midrule
solvent & \ce{MEA}, \ce{H2O}, \ce{CO2}, \ce{MEAH+}, \ce{MEACOO-}, \ce{H+}, \ce{HCO3-}, \ce{OH-}, \ce{CO3^2-} \\
vapor  & \ce{N2}, \ce{O2}, \ce{H2O}, \ce{CO2} \\
\bottomrule
\end{tabular}
\end{table}

\subsection{Unit-operation models}
\label{SI:unitoperations}
Columns are modeled as a series of well‑mixed segments, while all other units are treated as a single equilibrium or flash stage. All unit models share the following generic constraints: component balances \cref{eq:component_balance}, energy balances \cref{eq:energy_balance}, species summation \cref{eq:summation}, vapor-liquid equilibrium \cref{eq:VLE}, and enthalpy models for all scenarios $s\in\mathcal{S}$, segments $j\in\mathcal{J}$, and species $c\in\mathcal{C}$.

\begin{equation}
    \dot{L}_{s,j}\,x_{c,s,j} + \dot{V}_{s,j}\,y_{c,s,j}
    = \dot{L}_{s,j+1}\,x_{c,s,j+1} + \dot{V}_{s,j-1}\,y_{c,s,j-1}
    + \dot{F}^{\mathrm{L}}_{c,s,j} + \dot{F}^{\mathrm{V}}_{c,s,j}
    + \sum_{r\in\mathcal{R}}\nu_{c,r}\,\xi_{s,j,r},
    \quad\forall\, c\in\mathcal{C},\, s\in\mathcal{S},\, j\in\mathcal{J},
    \label{eq:component_balance}
\end{equation}
where $\dot{L}_{s,j}$ and $\dot{V}_{s,j}$ are the molar flow rates of liquid and vapor leaving segment $j$, $x_{c,s,j}$ and $y_{c,s,j}$ the liquid and vapor mole fractions of species $c$, $\dot{F}^{\mathrm{L}}_{c,s,j}$ and $\dot{F}^{\mathrm{V}}_{c,s,j}$ external liquid and vapor feed molar flow rates of species $c$ to segment $j$, $\nu_{c,r}$ the stoichiometric coefficient of species $c$ in reaction $r$, and $\xi_{s,j,r}$ the extent of reaction $r\in\mathcal{R}$ (\crefrange{eq:reaction1}{eq:reaction5}). For inert species, $\nu_{c,r}=0$.

The energy balance of each segment reads
\begin{equation}
    \begin{aligned}
        H^{\mathrm{L,out}}_{s,j} + H^{\mathrm{V,out}}_{s,j}
        = H^{\mathrm{L,in}}_{s,j} + &H^{\mathrm{V,in}}_{s,j} + \dot{Q}_{s,j}, 
        &\forall\, s\in\mathcal{S},\, j\in\mathcal{J},
    \end{aligned}
    \label{eq:energy_balance}
\end{equation}
where $H^{\mathrm{L,out}}_{s,j}$ and $H^{\mathrm{V,out}}_{s,j}$ are the total liquid and vapor enthalpies of streams leaving segment $j$, $H^{\mathrm{L,in}}_{s,j}$ and $H^{\mathrm{V,in}}_{s,j}$ those of streams entering segment $j$ (including any feed streams), and $\dot{Q}_{s,j}$ an external heat duty, nonzero only in units and stages with external heating or cooling.

The mole fractions in each phase sum to one:
\begin{equation}
    \sum_{c\in\mathcal{C}} x_{c,s,j}=1, \quad \sum_{c\in\mathcal{C}} y_{c,s,j}=1,
    \quad\forall\, s\in\mathcal{S},\, j\in\mathcal{J}.
    \label{eq:summation}
\end{equation}

For vapor-liquid equilibrium, we assume water to follow Raoult's law and \ce{CO2} to follow Henry's law:
\begin{equation}
    \begin{aligned}
        p\,y_{\ce{H2O},s,j} &= x_{\ce{H2O},s,j}\,p^{\mathrm{sat}}_{\ce{H2O}}(T_{s,j}), &\forall\, s\in\mathcal{S},\, j\in\mathcal{J},\\
        p^{*}_{\ce{CO2},s,j} &= \mathrm{He}_{\ce{CO2}}\,c_{\ce{CO2},s,j}, &\forall\, s\in\mathcal{S},\, j\in\mathcal{J}, 
    \end{aligned}
    \label{eq:VLE}
\end{equation}
where $p$ is the operating pressure, $p^{\mathrm{sat}}_{\ce{H2O}}(T_{s,j})$ the saturation pressure of water at temperature $T_{s,j}$, $p^{*}_{\ce{CO2},s,j}$ the equilibrium partial pressure of \ce{CO2}, $\mathrm{He}_{\ce{CO2}}$ Henry's constant for \ce{CO2}, and $c_{\ce{CO2},s,j}$ the molar concentration of dissolved \ce{CO2} in the liquid phase. Vapor-liquid equilibrium is enforced in all units except the columns by setting $p\,y_{\ce{CO2},s,j} = p^{*}_{\ce{CO2},s,j}$; in the columns, the approach-to-equilibrium model (\cref{eq:approach_to_equilibrium,eq:calibrationModel}) is used to account for rate limitations.

We assume ideal‑mixture enthalpies with pure‑component correlations for heat capacity, heat of vaporization, and constant heat of absorption collected by \citet{Akula.2021}, neglecting heat of mixing.

\subsubsection{Columns}
Column diameters are sized based on generalized pressure‑drop correlations for structured packings \citep{Tsai.2010}, limiting the superficial gas velocity to below \SI{80}{\percent} of flooding velocity in all operating scenarios \citep{Kister.2007}, assuming Mellapak 252Y. A minimum liquid load of \SI{0.2}{\cubic\meter\per\square\meter\per\hour} \citep{SulzerLtd.2024} enforces sufficient wetting. An intercooler is placed at mid‑height with the cooled‑segment liquid temperature constrained to at most \SI{50}{\celsius}. As a simplifying assumption, we assume absorption and desorption pressure levels to be fixed at \SI{110}{\kilo\pascal} and \SI{200}{\kilo\pascal}, respectively. The cooling duty of the intercooler is an operational variable and the maximum duty determines the sizing of the intercooler heat-exchanger.

For the costing of the columns, bare erected costs are estimated based on textbook correlations of \citet{Seider.2017}, assuming stainless steel equipment, escalated to 2023 with the Chemical Engineering Plant Cost Index. To estimate the column costs, we add an additional height required for water wash and sump based on the packing height assuming the same ratio as in a front-end engineering design study of natural gas carbon capture \citep{ElectricPowerResearchInstitute.2021} of 1.27 and 0.85 for the absorber and desorber, respectively. The cost of the structured packing material is estimated based on \citet{Tsai.2010}.

\subsubsection{Heat exchangers}
The cross heat exchanger is a four‑pass plate‑and‑frame unit modeled with the effectiveness–NTU method following \citet{akula.2019}. We assume a constant thermal conductivity of \SI{15.71}{\watt\per\meter\per\kelvin} in the temperature range of \SIrange{300}{400}{\kelvin} assuming AISI 304 stainless steel plates \citep{Bogaard.1985}, and a corrugation factor of 1.2. The plate sizing and gap are parametrized following \citet{Akula.2021}. 

The reboiler, condenser, lean‑amine cooler and intercooler are costed as shell‑and‑tube exchangers sized based on their respective maximum duties, the logarithmic mean temperature difference expressed using Chen's approximation \citep{Chen.1987}, and assuming constant heat-transfer coefficients (\Cref{tab:HX}) based on the detailed studies of \citet{Fisher.2007}. We assume saturated steam at \SI{150}{\psig} for regeneration and cooling water at \SI{90}{\fahrenheit} with a maximum outlet temperature of \SI{120}{\fahrenheit}. We assume an outlet temperature of the condenser of \SI{40}{\celsius}. The intercooler outlet temperature of the solvent is modeled as a variable constrained to under \SI{50}{\celsius}.  

\begin{table}[!t]\centering
\caption{Overall heat‑transfer coefficients $U$ based on \citet{Fisher.2007}.}
\label{tab:HX}
\begin{tabular}{ll}
\toprule
exchanger & $U$ [\si{\watt\per\square\meter\per\kelvin}] \\
\midrule
reboiler & 852 \\
condenser & 454 \\
lean‑amine cooler, intercooler & 795 \\
\bottomrule
\end{tabular}
\end{table}

\subsubsection{Parallelization of oversized equipment}
The large flue-gas mass flows in our case studies require heat-exchanger areas and pump flow rates that exceed the maximum sizes available in the cost correlations of \citet{Seider.2017}. To accommodate these large heat and volume flows, we assume multiple maximum-size units operating in parallel. To preserve the nonlinear programming formulation, we relax the number of units to a continuous variable. This relaxation introduces a limited modeling error: When we recalculate the costs assuming a discrete number of equipment for comparison, the bare erected costs of the heat exchangers and pumps are underestimated by \SI{1.4}{\percent} and \SI{1.2}{\percent}, respectively, in the base case of the coal case study.

\FloatBarrier
\subsection{Approach-to-equilibrium column model}
\label{SI:approachtoequilibrium}
The approach-to-equilibrium model approximates rate-limitations from chemical reactions and mass transfer. The form of the function $\eta(\cdot)$ \cref{eq:calibrationModel} is chosen such that $\eta \in [0,1]$. Further, $\eta$ is consistent across discretization granularities, and approaches 1 for $\Delta z \longrightarrow \infty$, and 0 for $\Delta z \longrightarrow 0$. 

The function $k(\cdot\,;\boldsymbol{\beta})$ of the approach-to-equilibrium model \cref{eq:approach_to_equilibrium} is evaluated using inputs normalized by their mean values $\mu_n$ and standard deviations $\sigma_n$: $\hat\zeta_{n,s,j}=(\zeta_{n,s,j}-\mu_n)/\sigma_n$. 

The coefficients $\beta_0$, $\beta_n$, and $\beta_{n,m}$ are fitted separately for the absorber and the desorber (\Cref{tab:fitting_parameters}) using calibration data simulated using the rate-based \textit{CCSI steady state MEA model} \citep{U.S.DepartmentofEnergy.2022} in Aspen Plus. The calibration data is generated by varying the inlet conditions of the liquid and vapor streams over the ranges in \Cref{tab:operatingrange}, which were selected in discussion with an industry partner. Liquid and vapor conditions are extracted using the Python package \textit{pywin32} at the column outlets and within the column. Simulation runs that did not converge were discarded prior to fitting. The coefficients $\boldsymbol{\beta}$ are then calibrated by minimizing the mean squared error of the predicted values of $\ln k(\cdot;\boldsymbol{\beta})$ using the Broyden–Fletcher–Goldfarb–Shanno algorithm as implemented in the Julia package \textit{Optim} \citep{KMogensen.2018}. 

The coal and gas case studies differ in flue-gas composition, so the absorber coefficients are fitted to calibration data with different \ce{CO2} concentration ranges (\Cref{tab:operatingrange}). 

\begin{figure*}[!t]
    \centering
    \includegraphics[width=5in]{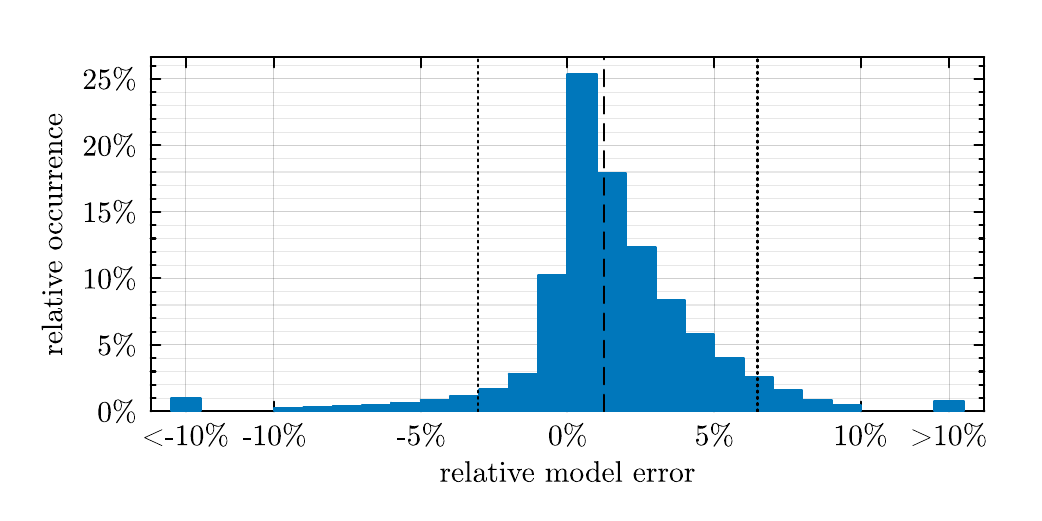}
    \caption{Prediction error of the approach-to-equilibrium column model's column segment outlet \ce{CO2} mole fraction compared to the simulated data from the CCSI steady state model \citep{U.S.DepartmentofEnergy.2022} for the case study treating coal power plant flue gases. The error is evaluated for over 360,000 column segment inlet conditions sampled from the range of operating conditions tabulated in \Cref{tab:operatingrange}. The dotted lines indicate the 5th and 95th percentile, while the dashed line indicates the average error.}
    \label{fig:predictionerror_coal}
\end{figure*}

\begin{figure*}[!ht]
    \centering
    \includegraphics[width=5in]{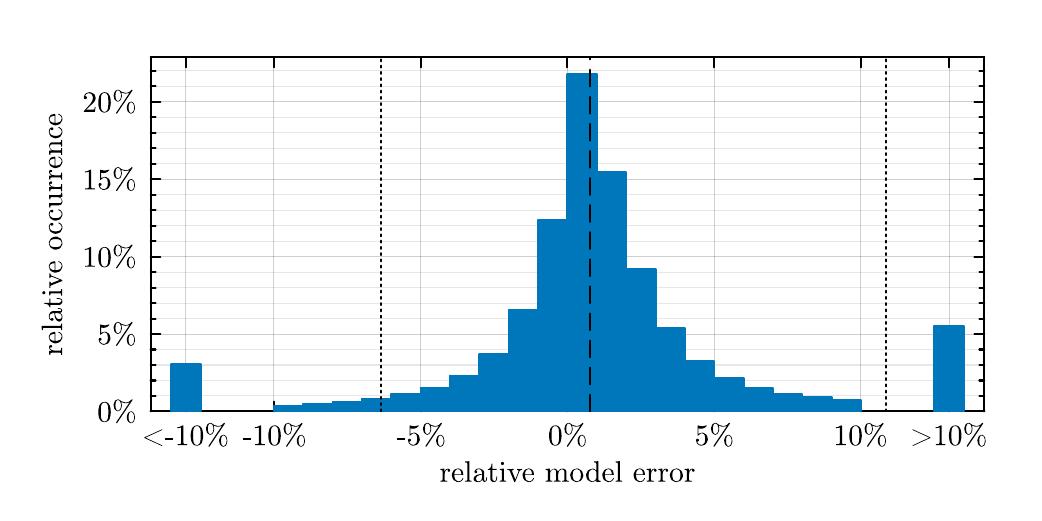}
    \caption{Prediction error of the approach-to-equilibrium column model's column segment outlet \ce{CO2} mole fraction compared to the simulated data from the CCSI steady state model \citep{U.S.DepartmentofEnergy.2022} for the case study treating natural gas combined cycle power plant flue gases. The error is evaluated for over 190,000 column segment inlet conditions sampled from the range of operating conditions tabulated in \Cref{tab:operatingrange}. The dotted lines indicate the 5th and 95th percentile, while the dashed line indicates the average error.}
    \label{fig:predictionerror_gas}
\end{figure*}

The relative prediction errors for the \ce{CO2} outlet mole fractions of an absorber stage for the calibated model are shown in \Crefrange{fig:predictionerror_coal}{fig:predictionerror_gas}.

\begin{table*}[!ht]\centering
\caption{Fitted coefficients, feature means, and standard deviations for the
absorber and desorber approach-to-equilibrium models \cref{eq:fittingFunction}.
Separate absorber models are fitted for coal and natural gas flue gas.}
\label{tab:fitting_parameters}
    \begin{tabular}{l rrr c rrr c rrr}
    \toprule
     & \multicolumn{7}{c}{absorber} && \multicolumn{3}{c}{desorber} \\
    \cmidrule{2-8}\cmidrule{10-12}
     & \multicolumn{3}{c}{coal} && \multicolumn{3}{c}{natural gas} && & & \\
    \cmidrule{2-4}\cmidrule{6-8}
    $n$ & $\beta_n$ & $\mu_n$ & $\sigma_n$ 
        && $\beta_n$ & $\mu_n$ & $\sigma_n$ 
        && $\beta_n$ & $\mu_n$ & $\sigma_n$ \\
    \midrule
    $0$
        & -2.714 & ---    & ---   
        && -2.926 & ---    & ---   
        && -2.686 & ---    & ---   \\
    $V/L$
        & -0.481 & 0.460  & 0.290 
        && -0.217 & 0.464  & 0.294 
        &&  0.795 & 0.100  & 0.063 \\
    $T_\text{in}$
        &  0.287 & 54.93  & 8.913 
        && -0.167 & 49.67  & 7.405 
        &&  1.084 & 109.65 & 8.548 \\
    $p$
        &  0.024 & 1.086  & 0.081 
        &&  0.018 & 1.086  & 0.081 
        && -0.533 & 1.776  & 0.083 \\
    $\alpha$
        &  0.313 & 0.339  & 0.125 
        &&  1.358 & 0.406  & 0.081 
        && -41.14 & 0.312  & 0.130 \\
    $x^{\ce{MEA},\mathrm{lean}}$
        &  0.152 & 0.302  & 0.005 
        &&  2.612 & 0.299  & 0.064 
        && -3.078 & 0.230  & 0.016 \\
    $u_v$
        & -0.198 & 2.119  & 0.779 
        && -0.192 & 2.048  & 0.746 
        && -0.015 & 0.638  & 0.526 \\
    \midrule
    $n,m$ & $\beta_{n,m}$ & $\mu_{n,m}$ & $\sigma_{n,m}$ 
          && $\beta_{n,m}$ & $\mu_{n,m}$ & $\sigma_{n,m}$ 
          && $\beta_{n,m}$ & $\mu_{n,m}$ & $\sigma_{n,m}$ \\
    \cmidrule{1-4}\cmidrule{6-8}\cmidrule{10-12}
    $V/L,\, u_v$
        &  0.044 & 0.961  & 0.707 
        &&  0.017 & 0.944  & 0.710 
        && -0.606 & 0.096  & 0.109 \\
    $T_\text{in},\, \alpha$
        & -0.335 & 18.23  & 6.355 
        &&  0.041 & 19.98  & 4.566 
        &&  16.70 & 33.12  & 11.83 \\
    $\alpha,\, x^{\ce{MEA},\mathrm{lean}}$
        & -0.821 & 0.103  & 0.038 
        && -3.303 & 0.121  & 0.037 
        &&  26.28 & 0.074  & 0.035 \\
    \bottomrule
    \end{tabular}
\end{table*}

\begin{table*}[!ht]
    \centering
    \caption{Operating ranges assumed for simulating calibration data used for fitting the approach-to-equilibrium model.}
    \begin{tabular}{lrr}
        \toprule
        parameter & unit & range \\
        \midrule
        input liquid-to-gas ratio & kg/kg & 1--5 \\
        input solvent temperature & \si{\degreeCelsius} & 40--70 \\
        absorber pressure & \si{\kilo\pascal} & 100--120 \\
        input solvent \ce{CO2} loading & -- & 0.1--0.5 \\
        input solvent composition & \si{\wtpercent} & 20--40 \\
        input vapor \ce{CO2} mole fraction (coal) & -- &  0.03--0.15\\
        input vapor \ce{CO2} mole fraction (natural gas) & -- &  0.01--0.05\\
        superficial gas velocity & \si{\meter\per\second} & 1--3 \\
        \bottomrule
    \end{tabular}
    \label{tab:operatingrange}
\end{table*}

\FloatBarrier

\subsection{Physical property correlations}
\label{SI:properties}
Physical property correlations applied in the carbon capture process model are listed in \Cref{tab:property_models}.
\begin{table}[!ht]
\centering
\caption{Physical‑property models and sources.}
    \begin{tabular}{ll}
    \toprule
    property & source \\
    \midrule
    \ce{H2O} vapor pressure & \citet{Perry.2008} \\
    \ce{CO2} Henry's constant & \citet{Hartono.2014} \\
    Liquid heat capacity  & \citet{Akula.2021} \\
    Liquid viscosity & \citet{Akula.2021} \\
    Thermal conductivity (\ce{MEA}) & \citet{DiGuilio.1992} \\
    Thermal conductivity (\ce{H2O}) & \citet{Ramires.1995} \\
    \bottomrule
    \end{tabular}
\label{tab:property_models}
\end{table}

As a simplifying assumption, liquid density is assumed constant at \SI{1000}{\kilo\gram\per\cubic\meter} and vapor density is calculated assuming ideal gases.

\FloatBarrier
\subsection{Economic parameters}
\label{SI:cost}
Bare erected costs of equipment are estimated via cost correlations by \citet{Seider.2017} listed in \ref{tab:costmodels}. Bare erected costs are converted to total plant cost by adding engineering, procurement, and  construction costs, as well as process and project contingencies, using cost factors. Total overnight cost further includes  preproduction, inventory, financing, and other costs, using factors from \citet{Cvetic.2023} following the National Energy Technology Laboratory capital cost estimation methodology \citep{Theis.2021}. The total overnight cost is annualized with the capital charge factor of \citet{Theis.2021}.

\begin{center}
\begin{longtable}{p{0.10\textwidth} p{0.06\textwidth} p{0.54\textwidth} p{0.20\textwidth}}
    \caption{bare erected cost correlations for equipment in the carbon capture process.}
    \label{tab:costmodels}
    \endfirsthead
    \endhead
    \toprule
    equipment & variable & cost calculation &  comments\\
    \midrule
    column shell & $d$, $h$ &   \textbf{column height} $h$ (\si{\meter}):\newline
                                $h = h^{\mathrm{pack}} (1+c^{\mathrm{add}})$\newline
                                \textbf{wall thickness} $t$ (in):\newline
                                $t = t^{\mathrm{base}} + t^{\mathrm{we}} + t^{\mathrm{corr}}$,\quad
                                $t^{\mathrm{base}} = 0.5$,\newline
                                with the allowance for corrosion, 
                                $t^{\mathrm{corr}} = 0.125$, and for wind and earthquakes $t^{\mathrm{we}}$:\newline
                                $t^{\mathrm{we}} = 0.22\,(39.37\, d+18)\,(39.37\,h)^{2} /   (\sigma_{\max}\,(39.37\,d)^{2})$,\quad$\sigma_{\max} = 20000\,$psi\newline
                                \textbf{shell weight} $w$ (\si{\kilo\gram}):\newline
                                $w = \rho^{s}\,\pi\,(d+t/\SI{39.37}{\inch\per\meter})\,(0.8\,d + h)\,t/\SI{39.37}{\inch\per\meter},$\newline
                                $\rho^{s} = \SI{8000}{\kilo\gram\per\cubic\meter}$\newline
                                \textbf{bare erected cost} $BEC^{\mathrm{shell}}$ (2013\,USD):\newline
                                $BEC^{\mathrm{shell}} = f^\mathrm{BM}\,\newline
                                \hspace*{5mm}\cdot\Big(f^{\mathrm{M}}\,\exp\!\big(10.545 - 0.467\ln \left(2.205\, w\right) \newline
                                \hspace*{5mm}+ 0.0548\,f^{\mathrm{M}}\,\left(\ln \left(2.205\, w \right)\right)^{2} \big)\newline
                                \hspace*{5mm}+ 341\,\left(3.28\,d\right)^{0.633}\,\left(3.28\,h\right)^{0.802}\Big)$, \newline
                                with $f^\mathrm{M} = 1.7$ (SS\,304L) and $f^\mathrm{BM} = 4.16$\newline
                                & $h^{\mathrm{pack}}$ is the packing height in \si{\meter}, $c^\mathrm{add}$ is 1.27 and 0.85 for absorber and desorber respectively \citep{ElectricPowerResearchInstitute.2021}, diameter $d$ and height $h$ in \si{\meter}, cost correlation based on \citet{Seider.2017} \\
    \addlinespace[6pt]
    packing & $d$, $h$ &        \textbf{bare erected cost} $BEC^{\mathrm{pack}}$                                            (2008\,USD):\newline
                                $BEC^{\mathrm{pack}} = f^\mathrm{BM}\, c^{\mathrm{v}}\,\pi/4\,d^{2}\,h$,\newline
                                with $c^{\mathrm{v}} = \SI{50}{\USD\per\cubic\foot}$, $f^\mathrm{BM} = 4.16$.
                                & diameter $d$ and height $h$ in \si{\meter}, specific cost from \citet{Tsai.2010} \\
    \addlinespace[6pt]
    pump & $\dot{v}$ &   \textbf{size factor} $SF$:\newline
                                $SF = \dot{v}\,H^{0.5}$,\newline
                                \textbf{pump efficiency} $\eta^{\mathrm{P}}$:\newline
                                $\eta^{\mathrm{P}} = -0.316 + 0.240\ln\dot{v} - 0.0120\,(\ln\dot{v})^{2}$\newline
                                \textbf{brake power} $P^{\mathrm{B}}$ (hp):\newline
                                $P^{\mathrm{B}} = \dot{v}\,H\,\rho\,/\,(33000\,\eta^{\mathrm{P}})$,\quad $\rho$ in lb\,gal$^{-1}$\newline
                                \textbf{motor efficiency} $\eta^{\mathrm{M}}$:\newline
                                $\eta^{\mathrm{M}} = 0.800 + 0.0319\ln P^{\mathrm{B}} - 0.00182\,(\ln P^{\mathrm{B}})^{2}$\newline
                                \textbf{motor power} $P^{\mathrm{M}}$ (hp):\newline
                                $P^{\mathrm{M}} = P^{\mathrm{B}}/\eta^{\mathrm{M}}$\newline
                                \textbf{bare erected cost} $BEC^{\mathrm{pump}}$ (2013\,USD):\newline
                                $BEC^{\mathrm{pump}} = f^\mathrm{BM}\,f^{M}\,f^{T}\,\newline
                                \hspace*{5mm}\cdot\exp\left(12.166 - 1.145\ln\mathrm{SF} + 0.0862\,(\ln\mathrm{SF})^{2}\right)\newline
                                \hspace*{5mm}+ f^\mathrm{BM}\,\exp\big(5.933 + 0.168\ln P^{\mathrm{M}} - 0.110\,(\ln P^{\mathrm{M}})^{2}\newline
                                \hspace*{5mm}+ 0.0714\,(\ln P^{\mathrm{M}})^{3} - 0.00638\,(\ln P^{\mathrm{M}})^{4}\big)$,\newline
                                with $f^{M} = 2.0$,\quad $f^{T} = 2.0$, $f^\mathrm{BM} = 3.3$
                                & centrifugal variable speed pumps, maximum volume flow $\dot{v}$ in \si{\gpm}, head $H$ in \si{\foot}, cost correlation based on \citet{Seider.2017} \\
    \addlinespace[6pt]
    cross heat exchanger &      $A$ &
                                \textbf{bare erected cost} $BEC^{\mathrm{hex}}$ (2013\,USD):\newline
                                $BEC^{\mathrm{hex}} = f^\mathrm{BM}\,10070\,A^{0.42}$,\quad $A$ in ft$^{2}$, with $f^\mathrm{BM} = 3.17$
                                & plate-and-frame, surface area $A$ in \si{\square\foot}, cost correlation based on \citet{Seider.2017} \\
    \addlinespace[6pt]
    shell-and-tube heat exchanger & $A$ &
                                \textbf{bare erected cost} $BEC^{\mathrm{s\&t}}$ (2013\,USD):\newline
                                $BEC^{\mathrm{s\&t}} = f^\mathrm{BM}\,\left((2.7 + (A/100)^{0.07})\right)\newline
                                \hspace*{5mm}\cdot\exp\left(11.419 - 0.9228\ln A + 0.09861\,(\ln A)^{2}\right)$,\newline 
                                with $f^\mathrm{BM} = 3.17$
                                & used for reboiler, condenser, lean-amine cooler, intercooler; surface area $A$ in \si{\square\foot}; cost correlation based on \citet{Seider.2017}, we apply an additional cost factor for the reboiler $f^{\mathrm{R}} = 1.35$ based on \citet{Woods.2007}\\
    \addlinespace[6pt]
    compressor & $P$ &   \textbf{bare erected cost} $BEC^{\mathrm{comp}}$ (2013\,USD):\newline
                                $BEC^{\mathrm{comp}} = f^{\mathrm{M}}\,\exp(9.155)\,P^{0.63}$,\newline
                                with $f^\mathrm{BM} = 2.15$, $f^{\mathrm{M}} = 2.5$ 
                                & multi-stage centrifugal, power $P$ in hp, cost correlation based on \citet{Seider.2017}\\
    \bottomrule
\end{longtable}
\end{center}

Operating costs comprise 1) utility costs (\SI{150}{\psig} steam, electricity for pumps and compressor, and \SI{90}{\fahrenheit} cooling water, priced following \citet{Seider.2017}), 2) fixed costs for labor, maintenance, taxes, and insurance \citep{Cvetic.2023}, and 3) variable costs for waste treatment, makeups, and chemicals scaled with the amount of \ce{CO2} captured \citep{Cvetic.2023}. The \ce{MEA} makeup cost is based on the degradation rates of \citet{Moser.2020} and prices of \citet{Haslbeck.2010}, escalated to 2023 using the producer price index for industrial chemicals \citep{U.S.BureauofLaborStatistics.2026}. The main economic parameters are presented in \Cref{tab:economics}.

\begin{table}[!t]\centering
\caption{Key economic parameters assumed in the carbon capture process model}
    \begin{tabular}{p{0.25\textwidth}p{0.45\textwidth}p{0.2\textwidth}}
    \toprule
    parameter & value & source \\
    \midrule
    capital charge factor & 0.0773 & \citet{Theis.2021} \\
    engineering, procurement, and construction & \SI{20}{\percent} of bare erected cost & \citet{Cvetic.2023} \\
    process / project contingency & \SI{18}{\percent} / \SI{27.6}{\percent} of bare erected cost & \citet{Cvetic.2023} \\
    retrofit factor & 1.05 & \citet{Cvetic.2023} \\
    steam (\SI{150}{\psig}) & \SI{15.3}{USD\per\tonne} & \citet{Seider.2017} \\
    electricity & \SI{70}{USD\per\mega\watt\hour} & \citet{Seider.2017} \\
    cooling water (\SI{90}{\fahrenheit}) & \SI{0.027}{USD\per\cubic\meter} & \citet{Seider.2017} \\
    \ce{MEA} (2007) & \SI{2480}{USD\per\tonne} & \citet{Haslbeck.2010} \\
    labor & \SI{38.5}{USD\per\hour}, 1.3 operators, \SI{55}{\percent} labor benefits overheads, and other burdens & \citet{Cvetic.2023} \\
    maintenance labor / material & \SI{0.76}{\percent} / \SI{1.14}{\percent} of total plant cost & \citet{Cvetic.2023} \\
    taxes and insurance & \SI{2}{\percent} of total plant cost & \citet{Cvetic.2023} \\
    \bottomrule
    \end{tabular}
\label{tab:economics}
\end{table}

\FloatBarrier
\section{Validation cases for flue gases with low-\ce{CO2} concentration}
The inlet conditions of flue gas and lean solvent during measurement campaigns at the Technology Centre Mongstad are collected from the literature for the validation of the approach-to-equilibrium model for the \gls{ngcc} case study (\Cref{tab:validation_gas_fluegas,tab:validation_gas_solvent}).

\begin{table*}[!ht]
    \centering
    \caption{Flue gas inlet conditions for the Technology Centre Mongstad validation cases.}
    \label{tab:validation_gas_fluegas}
        \begin{threeparttable}
        \small
        \begin{tabular}{l l l l l l l l p{4cm}}
            \toprule
            case
                & $\dot{V}$
                & $y_{\ce{CO2}}$
                & $y_{\ce{H2O}}$
                & $y_{\ce{O2}}$
                & $y_{\ce{N2}}$
                & $T$
                & $p$
                & source \\
            & [\SI{}{\meter\cubed\per\hour}] & [--] & [--] & [--] & [--]
            & [\degree C] & [kPa] & \\
            \midrule
            1 & 46970
               & 0.0370 & 0.0295 & 0.1360 & 0.7975
               & 25.0 & 106.3
               & \citep{Hamborg.2014, Oi.2021} \\
            2 & 59430
               & 0.0370 & 0.0370 & 0.1460 & 0.7800
               & 29.8 & 101.0
               & \citep{Faramarzi.2017} \\
            3 & 46900
               & 0.0349 & 0.0370\tnote{a} & 0.1460\tnote{a} & 0.7881
               & 28.2 & 101.0\tnote{a}
               & \citep{Gjernes.2017}\\
            4 & 46900
               & 0.0349 & 0.0370\tnote{a} & 0.1460\tnote{a} & 0.7881
               & 29.0 & 101.0\tnote{a}
               & \citep{Gjernes.2017}  \\
            5 & 59400
               & 0.0349 & 0.0370\tnote{a} & 0.1460\tnote{a} & 0.7881
               & 30.4 & 101.0\tnote{a}
               & \citep{Gjernes.2017}  \\
            6 & 59600
               & 0.0349 & 0.0370\tnote{a} & 0.1460\tnote{a} & 0.7881
               & 31.7 & 101.0\tnote{a}
               & \citep{Gjernes.2017}  \\
            \bottomrule
        \end{tabular}
        \begin{tablenotes}
            \footnotesize
            \item[a] Value assumed equal to case 2.
        \end{tablenotes}
    \end{threeparttable}
\end{table*}

\begin{table*}[!ht]
    \centering
    \caption{Lean solvent inlet conditions for Technology Centre Mongstad validation cases.}
    \label{tab:validation_gas_solvent}
    \begin{threeparttable}
        \small
        \begin{tabular}{l l l l l l l r}
            \toprule
            case
                & $\dot{m}$
                & $\dot{m}_{\ce{MEA}}/\dot{m}_{\ce{H2O}}$
                & $\dot{n}_{\ce{CO2}}/\dot{n}_{\ce{MEA}}$
                & $T$
                & source\\
            & [kg/h] & [mol/mol] & [mol/mol]
            & [$^\circ$C] & & \\
            \midrule
            1 & 54900 & 0.30 & 0.23 & 36.5
               & \citep{Hamborg.2014} \\
            2 & 57434 & 0.31 & 0.20 & 37.0
               & \citep{Faramarzi.2017} \\
            3 & 59600 & 0.31 & 0.29 & 37.0\tnote{a}
               & \citep{Gjernes.2017}  \\
            4 & 52500 & 0.31 & 0.25 & 37.0\tnote{a}
               & \citep{Gjernes.2017}  \\
            5 & 59200 & 0.31 & 0.26 & 37.0\tnote{a}
               & \citep{Gjernes.2017}  \\
            6 & 67700 & 0.31 & 0.22 & 37.0\tnote{a}
               & \citep{Gjernes.2017}  \\
            \bottomrule
        \end{tabular}
        \begin{tablenotes}
            \footnotesize
            \item[a] Value assumed equal to case 2.
        \end{tablenotes}
    \end{threeparttable}
\end{table*}

\FloatBarrier
\section{Flue gas characteristics of the coal power plant}
\label{SI:part_load}
The flue gas mass flows and \ce{CO2} concentrations for the coal power plant are extracted from \citet{Roeder.2014} (\Cref{tab:coal_flue_gas_data}). We use the flue gas mass flows and \ce{CO2} concentrations for the \gls{ngcc} power plant reported by \citet{Alcaraz-Calderon.2019}.

\begin{table*}[!ht]
\centering
\caption{Coal power plant flue gas characteristics at different load levels \citep{Roeder.2014}}
\label{tab:coal_flue_gas_data}
\begin{tabular}{rll}
\toprule
load level & flue gas flow (kg/s) & \ce{CO2} content (\SI{}{\molpercent}) \\
\midrule
\SI{100}{\percent} & 564.6 & 13.7 \\
\SI{90}{\percent} & 520.46 & 13.4 \\
\SI{80}{\percent} & 474.4 & 13.0 \\
\SI{70}{\percent} & 427.0 & 12.7 \\
\SI{60}{\percent} & 383.4 & 12.1 \\
\SI{50}{\percent} & 343.0 & 11.2 \\
\SI{40}{\percent} & 298.5 & 10.3 \\
\bottomrule
\end{tabular}
\end{table*}

\newpage
\section{Results for natural gas combined cycle with intermediate load profile}

\begin{figure}[!ht]
    \centering
    \begin{subfigure}[t]{0.48\textwidth}
        \centering
        \includegraphics[width=\linewidth]{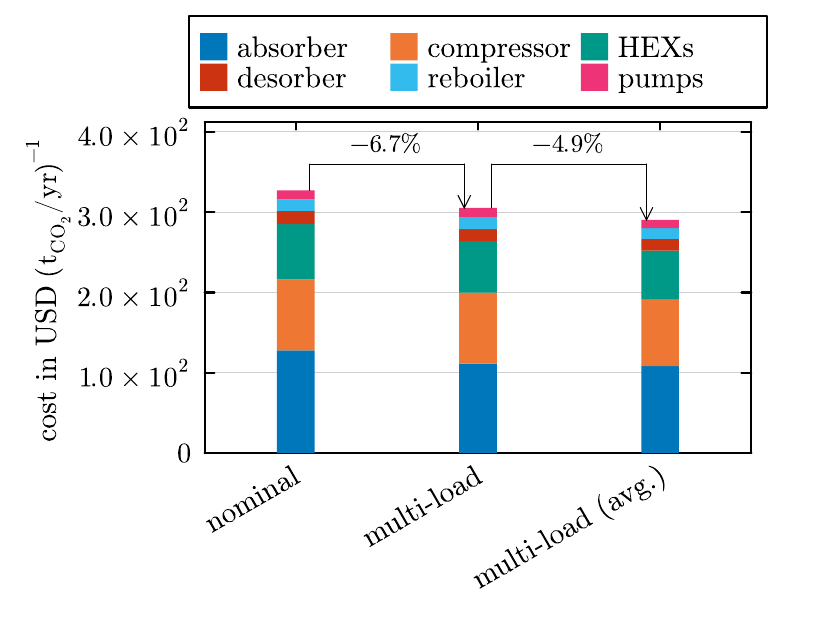}
        \caption{Total plant cost and main contributors. The designs of the two nominal cases are identical.}
        \label{fig:capex_gas_intermediate}
    \end{subfigure}
    \hfill
    \begin{subfigure}[t]{0.48\textwidth}
        \centering
        \includegraphics[width=\linewidth]{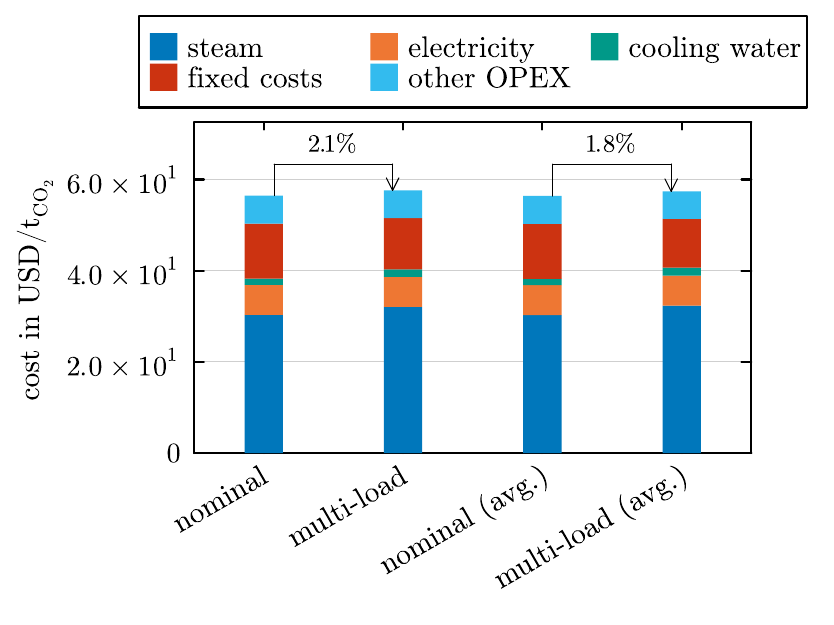}
        \caption{Annual operating cost and main contributors.}
        \label{fig:opex_gas_intermediate}
    \end{subfigure}
    
    \vspace{1.em}
    
    \begin{subfigure}[t]{0.48\textwidth}
        \centering
        \includegraphics[width=\linewidth]{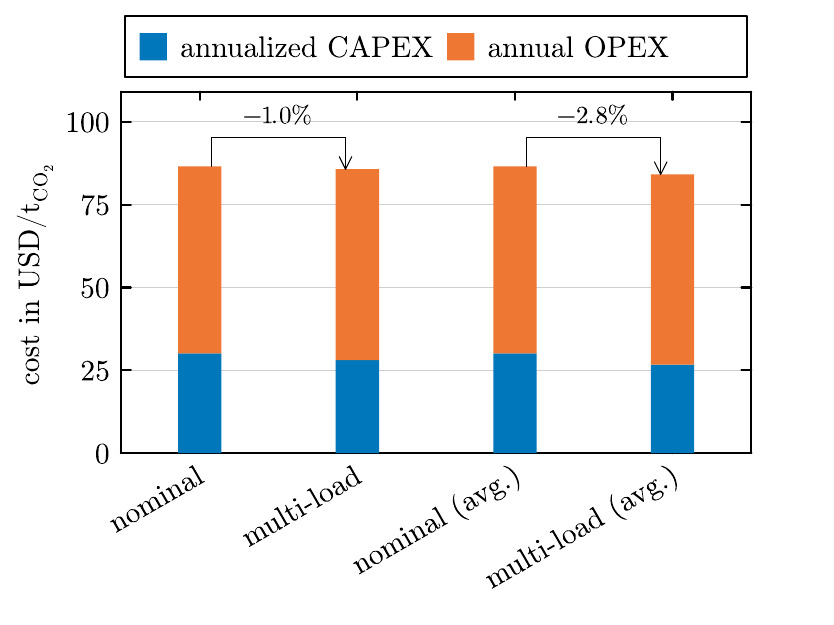}
        \caption{Total annualized cost of carbon capture.}
        \label{fig:TAC_gas_intermediate}
    \end{subfigure}
    \vspace{1.5em}
    \caption{Cost overview of capital and operating costs for the nominal and multi-load designs in the natural gas power plant with intermediate load (\Cref{fig:NGCC_medium}), assuming \SI{95}{\percent} capture rate in every scenario or on average over the full year (avg.). The capacity factor in the intermediate-load case is \SI{50.9}{\percent}. For comparison, in the high-load case, the capacity factor is \SI{59.7}{\percent}.}
    \label{fig:cost_gas_intermediate}
\end{figure}
\FloatBarrier

\newpage
\section{Sensitivity of cost savings}
\begin{figure*}[!ht]
    \centering
    \begin{subfigure}{0.72\textwidth}
        \centering
        \includegraphics[width=\textwidth,trim={0 1cm 0 0},clip]{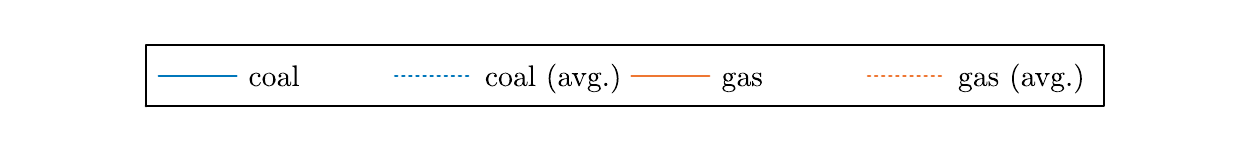}
    \end{subfigure}

    \begin{subfigure}[t]{0.48\textwidth}
        \centering
        \includegraphics[width=\textwidth,trim={0 0.5cm 0 0},clip]{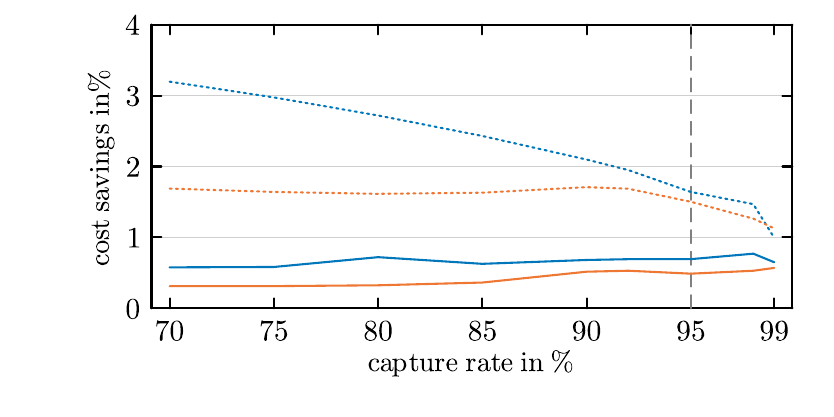}
        \caption{Sensitivity of cost savings to variations in capture rate target.}
        \label{fig:sensitivity_capture_rate_costsaving}
    \end{subfigure}
    \hfill
    \begin{subfigure}[t]{0.48\textwidth}
        \centering
        \includegraphics[width=\textwidth,trim={0 0.5cm 0 0},clip]{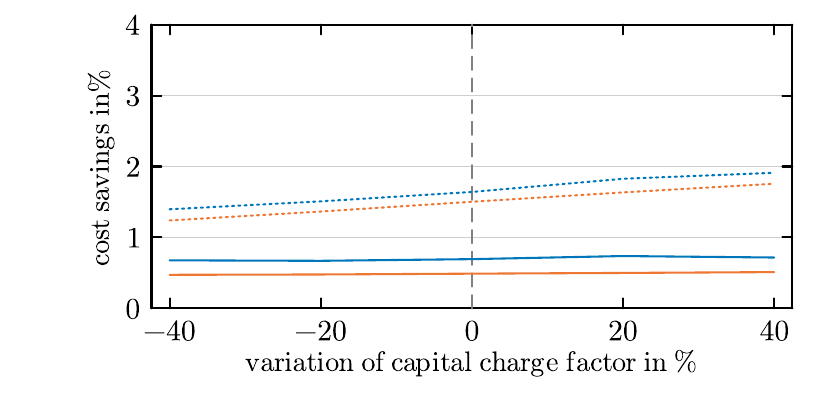}
        \caption{Sensitivity of cost savings to variations of capital charge factor.}
        \label{fig:sensitivity_CCF_costsaving}
    \end{subfigure}
    
    \vspace{0.3em}
    
    \begin{subfigure}[t]{0.48\textwidth}
        \centering
        \includegraphics[width=\textwidth,trim={0 0.5cm 0 0},clip]{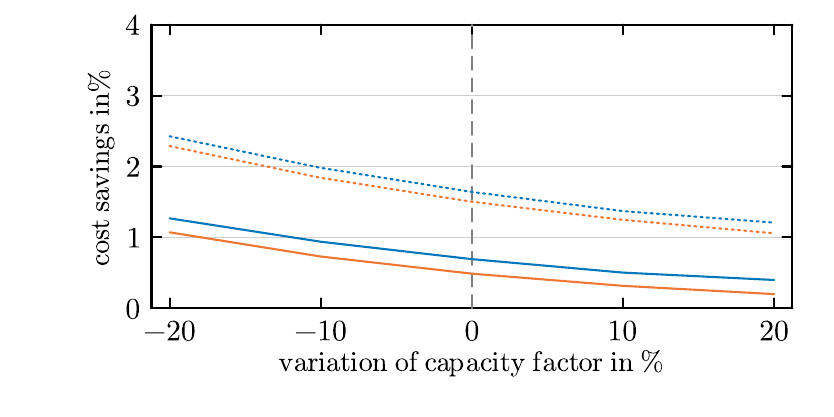}
        \caption{Sensitivity of cost savings to variations in capacity factor.}
        \label{fig:sensitivity_flh_costsaving}
    \end{subfigure}
    \hfill
    \begin{subfigure}[t]{0.48\textwidth}
        \centering
        \includegraphics[width=\textwidth,trim={0 0.5cm 0 0},clip]{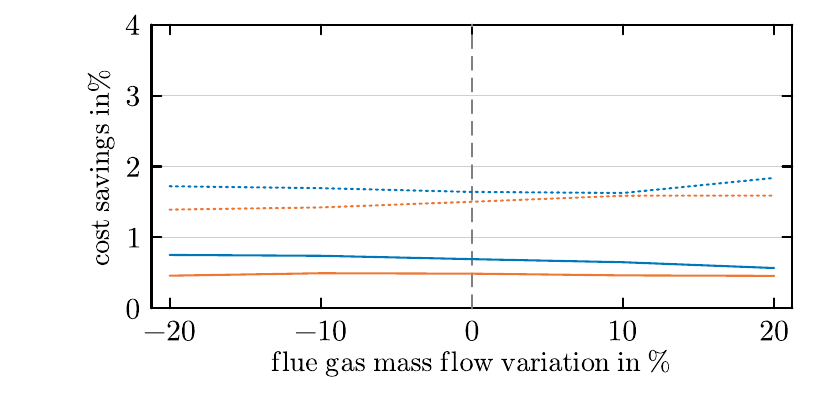}
        \caption{Sensitivity of cost savings to variations in flue gas mass flow.}
        \label{fig:sensitivity_massflow_costsaving}
    \end{subfigure}
    
    \caption{Sensitivity of the relative cost savings of the multi-load design in comparison to the nominal design for both coal and natural gas combined cycle case studies. Cost savings are presented for cases with a fixed capture rate per load scenario and for cases with an average capture rate (avg.) across scenarios}
    \label{fig:sensitivity_costsaving}
\end{figure*}
\newpage

\section{Capture rate per load scenario for average capture rate targets across scenarios}
\begin{figure*}[!h]
    \centering
    \begin{subfigure}{0.6\textwidth}
        \centering
        \includegraphics[width=1\textwidth,trim={0 0 0 0},clip]{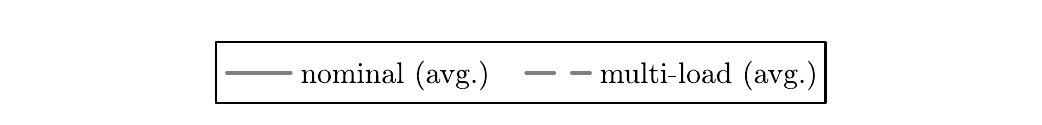}
    \end{subfigure}

    \begin{subfigure}[t]{0.48\textwidth}
        \centering
        \includegraphics[width=\textwidth,trim={0 0.5cm 0 0},clip]{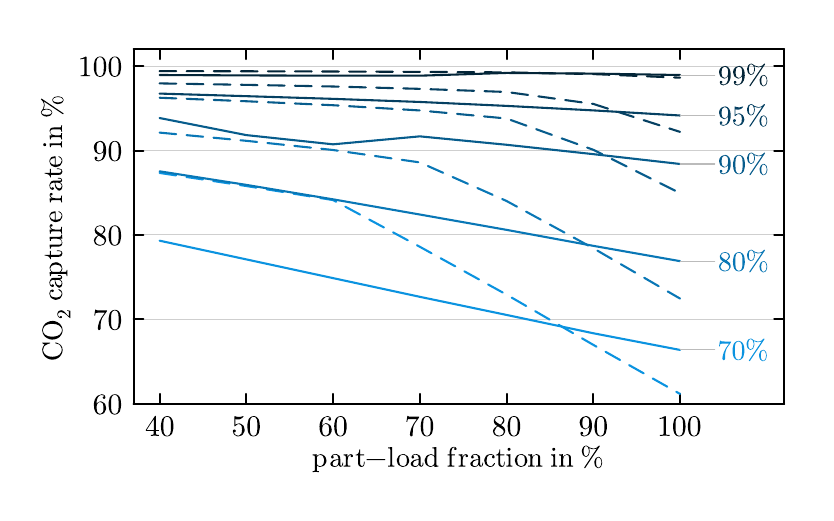}
        \caption{Capture rate depending on the part-load scenario for the coal power plant.}
        \label{fig:capture_rate_vs_part_load_coal}
    \end{subfigure}
    \hfill
    \begin{subfigure}[t]{0.48\textwidth}
        \centering
        \includegraphics[width=\textwidth,trim={0 0.5cm 0 0},clip]{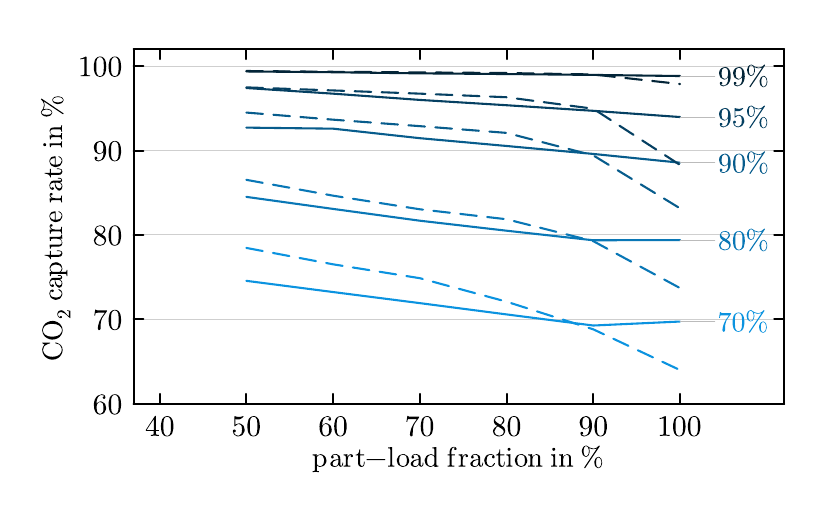}
        \caption{Capture rate depending on the part-load scenario for the natural gas combined cycle power plant.}
        \label{fig:capture_rate_vs_part_load_gas}
    \end{subfigure}
    \caption{Capture rate depending on the part-load scenario for coal and natural gas combined cycle case studies. The achieved capture rate per part-load scenario is shown for average capture targets ranging from \SI{70}{\percent} to \SI{99}{\percent} (see annotations in the figures) corresponding to the range of capture rates in the sensitivity analysis (\Cref{fig:sensitivity_capture_rate}).}
    \label{fig:capture_rate_vs_part_load}
\end{figure*}

\end{document}

%% file: preamble.tex
\usepackage[usenames,dvipsnames]{xcolor}
\usepackage{listings}
\usepackage{graphicx}
\usepackage{tikz}
\usetikzlibrary{positioning,backgrounds,fit,arrows.meta,calc}
\usepackage{pgfplots}
\pgfplotsset{compat=1.18}
\usepackage{amsmath,amssymb,amsthm,amsfonts}

\definecolor{mitred}{HTML}{750014}      
\definecolor{mitblue}{HTML}{1966FF}     
\definecolor{mitgreen}{HTML}{00AD00}    
\definecolor{mitdarkgray}{HTML}{40464C} 
\usepackage[linktocpage=true,colorlinks=true,linkcolor=blue,citecolor=blue,urlcolor=blue]{hyperref}
\usepackage[letterpaper,margin=1in]{geometry}
\usepackage{authblk}   

\makeatletter
\def\thm@space@setup{%
  \thm@preskip=0.6\baselineskip plus 2pt minus 2pt
  \thm@postskip=\thm@preskip}
\makeatother

\theoremstyle{definition}

\newtheorem{remark}{Remark}

\definecolor{jlkw}{HTML}{7B3294}      
\definecolor{jltype}{HTML}{1F78B4}    
\definecolor{jlmacro}{HTML}{E31A1C}   
\definecolor{jlstr}{HTML}{A6761D}     
\definecolor{jlcmt}{HTML}{636363}     
\definecolor{jlnum}{HTML}{1B7837}     
\definecolor{jlbg}{HTML}{F7F7F7}      

\lstdefinelanguage{Julia}{%
  morekeywords={abstract,break,case,catch,const,continue,do,else,elseif,%
    end,export,false,for,function,global,if,import,in,let,local,%
    macro,module,mutable,nothing,quote,return,struct,true,try,%
    type,typealias,using,where,while},
  morekeywords=[2]{Int,Int64,Int32,Float64,Float32,Bool,String,Symbol,%
    Vector,Matrix,Array,Dict,Tuple,NamedTuple,Any,Nothing,Type,%
    AbstractFloat,AbstractArray,AbstractVector,AbstractMatrix},
  morekeywords=[3]{@time,@views,@inbounds,@simd,@kernel,@index,@sprintf,%
    @info,@eval,@constraint,@variable,@objective},
  sensitive=true, alsoother={\$},
  morecomment=[l]\#, morecomment=[n]{\#=}{=\#},
  morestring=[s]{"}{"}, morestring=[m]{'}{'},
}[keywords,comments,strings]

\providecommand{\keywords}[1]{%
  \small\textbf{\textit{Keywords---}} #1%
}

\usepackage{algorithm}
\usepackage{algpseudocode}
\usepackage{parskip}
\usepackage[font=small,labelfont=bf]{caption}
\usepackage{subcaption}
\usepackage{booktabs}
\usepackage{multirow}
\usepackage{enumitem}
\setlist[itemize]{leftmargin=1.5em}
\setlist[enumerate]{leftmargin=1.8em}
\usepackage{mathtools}
\allowdisplaybreaks

\usepackage[style=numeric,backend=biber,natbib=true,doi=true,url=false,eprint=true,isbn=false,maxbibnames=6, minbibnames=6,sorting=none, giveninits=true]{biblatex}
\DeclareNameAlias{default}{family-given}
\usepackage[acronym,nomain,nonumberlist]{glossaries}
\glsdisablehyper   
\newacronym{nlp}{NLP}{nonlinear program}
\newacronym{gpu}{GPU}{graphics processing unit}
\newacronym{kkt}{KKT}{Karush--Kuhn--Tucker}
\newacronym{simd}{SIMD}{single instruction, multiple data}
\newacronym{ams}{AMS}{algebraic modeling system}
\newacronym{mea}{MEA}{monoethanolamine}
\newacronym{ngcc}{NGCC}{natural gas combined cycle}

\usepackage[capitalise,noabbrev]{cleveref}
\crefformat{equation}{(#2#1#3)}
\Crefformat{equation}{Equation~(#2#1#3)}
\crefmultiformat{equation}{(#2#1#3)}%
{ and~(#2#1#3)}{, (#2#1#3)}{ and~(#2#1#3)}
\crefrangeformat{equation}{(#3#1#4) to~(#5#2#6)}
    
\crefname{assumption}{assumption}{assumptions}
\Crefname{assumption}{Assumption}{Assumptions}
\crefname{algorithm}{algorithm}{algorithms}
\Crefname{algorithm}{Algorithm}{Algorithms}
\crefname{listing}{listing}{listings}
\crefname{section}{Section}{Sections}
\Crefname{section}{Section}{Sections}
\crefname{appendix}{Appendix}{Appendices}
\Crefname{appendix}{Appendix}{Appendices}

\usepackage{placeins}
\usepackage[version=4]{mhchem}
\usepackage{chemmacros}

\usepackage{siunitx}
\usepackage{gensymb} 
\DeclareSIPrefix\Mio{Mio.\,}{6}
\DeclareSIUnit{\volpercent}{vol.\%}
\DeclareSIUnit{\molpercent}{mol\%}
\DeclareSIUnit{\wtpercent}{wt.\%}
\DeclareSIUnit{\kgGHG}{kg\,\ce{CO2}-eq.}
\DeclareSIUnit{\gGHG}{g\,\ce{CO2}-eq.}
\DeclareSIUnit{\tGHG}{t\,\ce{CO2}-eq.}
\DeclareSIUnit{\TWh}{TWh}
\DeclareSIUnit\kWh{kWh}
\DeclareSIUnit\euro{EUR}
\DeclareSIUnit\MWh{MWh}
\DeclareSIUnit\year{a}
\DeclareSIUnit\inch{in}
\DeclareSIUnit\foot{ft}
\DeclareSIUnit\lb{lb}
\DeclareSIUnit\bar{bar}
\DeclareSIUnit\psig{psig}
\DeclareSIUnit\gpm{gpm}
\DeclareSIUnit\fahrenheit{\degree{}F}
\DeclareSIUnit\USD{USD}

\usepackage{longtable}
\usepackage{threeparttable}

\usepackage{setspace}